\documentclass[a4paper, 12pt]{article}

\usepackage{multicol}
\usepackage{ae}
\usepackage{aecompl}
\usepackage[cm]{aeguill}
\usepackage[T1]{fontenc}
\usepackage[latin1]{inputenc}
\usepackage[english]{babel}
\usepackage{fancyhdr}
\usepackage{amsmath,amssymb,amsthm}
\usepackage[all]{xy}
\usepackage{pstricks,pst-node}
\usepackage{geometry}
\renewcommand\c[1]{\mathcal{#1}}

\title{The universal cover of an algebra without double bypass}
\author{Patrick Le Meur}
\date{ }

\newtheorem{theorem}{Theorem}[section]
\newtheorem{proposition}[theorem]{Proposition}
\newtheorem{lemma}[theorem]{Lemma}
\newtheorem{corollary}[theorem]{Corollary}
\newtheorem{remark}{Remark}[]
\newtheorem{definition}[theorem]{Definition}
\newtheorem{THM}{Theorem}[]
\newtheorem{THMM}{Theorem}[]

\newtheorem{example}{Example}[]

\begin{document}
\maketitle{}
\abstract{
  Let $A$ be a basic and connected finite dimensional algebra over a field
  $k$ of characteristic zero. We show that if the quiver of $A$ has no
  double bypass then the fundamental group (as defined in
  \cite{martinezvilla_delapena}) of any presentation of $A$
  by quiver and relations is the quotient of the
  fundamental group of a privileged presentation of $A$.
  Then we show that the Galois covering of $A$ associated with this
  privileged presentation 
satisfies a universal property with respect to the connected Galois
coverings of $A$ in a similar fashion to 
   the universal cover of a topological space.
}
\section*{Introduction}
In this text, $k$ will be an algebraically closed field.
Let $A$ be a finite dimensional algebra over $k$. In order to
study left $A$-modules we may assume that $A$ is basic and connected,
where basic means that $A$ is the direct sum of pairwise non
isomorphic indecomposable projective left $A$-modules.
For such an algebra, the study of the
Galois coverings of $A$ gives some information on the representation
theory of $A$ (see \cite{cibils_marcos}, \cite{green} and \cite{martinezvilla_delapena}) and is a
particular case of the covering techniques introduced in
\cite{bongartz_gabriel}, \cite{gabriel} and \cite{riedtmann}.
In order to manipulate coverings of $A$ we will always
consider, unless otherwise stated,
$A$ as a $k$-category with set of
objects a complete set $\{e_i\}_i$ of primitive pairwise
orthogonal idempotents and with morphisms space $e_i\to e_j$ the
vector space $e_jAe_i$. The covering techniques have led to the
definition (see \cite{green} and \cite{martinezvilla_delapena})
of a fundamental group associated with any presentation of
$A$ by quiver and admissible relations, and which satisfies
many topological flavoured properties (see \cite{assem_delapena},
\cite{green} and
\cite{martinezvilla_delapena}).
This construction and its associated properties depend on the choice
of a presentation of $A$. In particular, one can find algebras for
which there exist different presentations giving rise to non
isomorphic fundamental groups.
In this text we compare the fundamental groups
of the presentations of $A$ as defined in
\cite{martinezvilla_delapena}, and we study the coverings of $A$
with the following question in mind: does $A$ have a universal Galois
covering? i.e. does $A$ admit a Galois covering which is
factorised by any other Galois covering?
This question has been successfully treated when $A$ is
representation-finite (see \cite{bongartz_gabriel} and \cite{gabriel}).
The present study will involve quivers ``without double bypass''. In simple
terms, a quiver without double
bypass is a quiver which has no distinct parallel
arrows, no oriented cycle and has no subquiver of the following form
\psset{unit=2pt}
\begin{pspicture}(32,12)
\rput(0,0){\rnode{A}{}}
\rput(8,6){\rnode{B}{}}
\rput(24,6){\rnode{C}{}}
\rput(32,0){\rnode{D}{}}
\ncline{->}{A}{D}
\ncline{->}{B}{C}
\ncline[linestyle=dotted]{->}{A}{B}
\ncline[linestyle=dotted]{->}{C}{D}
\ncarc[arcangleA=60,arcangleB=60,linestyle=dotted]{->}{B}{C}
\end{pspicture}
where continued (resp. dotted) arrows represent arrows (resp. oriented
paths) in the quiver.
Assuming that $k$ is a characteristic zero
field and that the ordinary quiver $Q$ of $A$ has no double bypass, we
 prove the following result announced in \cite{lemeur}:
\begin{THM}
  \label{A}
  Assuming the above conditions, there exists a presentation $kQ/I_0\simeq A$ by
  quiver and relations such that
  for any other presentation $kQ/I\simeq
  A$ the identity map on the walks in $Q$ induces a surjective group morphism $\pi_1(Q,I_0)\twoheadrightarrow \pi_1(Q,I)$.
\end{THM}
\noindent{The} proof of the above theorem allows us to recover the following fact
proven in \cite{bardzell_marcos}: if $A$ is a basic
triangular connected and constricted finite dimensional $k$-algebra,
then different presentations of $A$ give rise to isomorphic
fundamental groups. Under the hypotheses made before stating Theorem \ref{A} and with the
same notations, if $k\tilde{Q}/\tilde{I}_0\xrightarrow{F_0} kQ/I_0$ is the Galois
covering with group $\pi_1(Q,I_0)$ induced by the universal Galois
covering of $(Q,I_0)$ (see \cite{martinezvilla_delapena}), we show the
following result.
\begin{THM}
  \label{B}
  For any connected Galois covering $F\colon \c C'\to A$ with group $G$ there
  exist an isomorphism $kQ/I_0\xrightarrow{\sim} A$, a Galois
  covering $p\colon k\tilde{Q}/\tilde{I}_0\to \c C'$ with group a normal
  subgroup $N$ of $\pi_1(Q,I_0)$ and a commutative diagram:\\
  \null\hfill\xymatrix{
    k\tilde{Q}/\tilde{I}_0 \ar@{->}[r]^p \ar@{->}[d]_{F_0} &  \c C' \ar@{->}[d]^F \\
    kQ/I_0 \ar@{->}[r]^{\sim} & A
  }\hfill\null\\
  together with an exact sequence of groups: $
  1\to N\to \pi_1(Q,I_0)\to G\to 1 $
\end{THM}
Hence Theorem \ref{B} partially answers the question concerning the
existence of a universal Galois covering.
The text is organised as follows: in Section $1$ we define
the notions we will use, in Section $2$ we prove
Theorem \ref{A}, in Section $3$ we give useful results
on covering functors, these results will be used in the proof of Theorem
\ref{B} to which Section $4$ is devoted. Section $2$ gives the
proofs of all the results that were announced by the author in
\cite{lemeur}. This text is part of the
author's thesis made at Université Montpellier 2 under the supervision
of Claude Cibils.
\section{Basic definitions}
\label{prel}
\textbf{$k$-categories, covering functors, Galois coverings}\\
A \textbf{$k$-category} is a
category $\c C$ such that the objects class  $\c C_0$ of $\c C$ is a
non empty set and such that each set $_y\c C_x$ of morphisms $x\to y$ of $\c C$ is a
$k$-vector space with $k$-bilinear composition. Let $\c C $ be a
$k$-category.
We will say that $\c C$ is
\textbf{locally bounded} if the following  properties are
satisfied:\\
\indent{a)} distinct objects are not isomorphic,\\
\indent{b)} for each $x\in \c C_0$, the $k$-algebra $_x\c C_x$ is local,\\
\indent{c)} $\oplus_{y\in \c
    C_0}\ _y\c C_x$ is finite dimensional for any $x\in\c C_0$,\\
\indent{d)}  $\oplus_{x\in \c
    C_0}\ _y\c C_x$ is finite dimensional for any $y\in\c C_0$.\\
Unless otherwise stated, all the $k$-categories we will introduce
will be locally bounded. As an example, let $A$ be a basic finite
dimensional $k$-algebra. If
$1=\sum_{i=1}^ne_i$ is a decomposition of the unit into a sum of
primitive orthogonal idempotents, then $A=\oplus_{i,j}e_jAe_i$
and $A$ is a locally bounded $k$-category with set of objects $\{e_1,\ldots,e_n\}$ and
with morphisms space $e_i\to e_j$ equal to $e_jAe_i$.
We will say that the $k$-category
$\c C$ is \textbf{connected} if for any $x,y\in\c C_0$ there exists a
sequence $x_0=x,\ldots,x_n=y$ in $\c C_0$ such that $_{x_i}\c
C_{x_{i+1}}\neq 0$ or $_{x_{i+1}}\c C_{x_i}\neq 0$ for any $i$. Recall
that an \textbf{ideal} $I$ of $\c C$ is the data of vector subspaces $_yI_x\subseteq\
_y\c C_x$ for each
$x,y\in\c C_0$, such that the composition of a morphism in $I$ with
any morphism in $\c C$ lies in $I$.
The \textbf{radical} (see \cite{bongartz_gabriel}) of $\c C$ is the ideal $\c R\c C$ of
$\c C$ such that
$_y\c R\c C_x$ is the set of non invertible morphisms $x\to y$ for any
$x,y\in\c C_0$. If $n\geqslant 2$
we set $\c R^n\c C=\left(\c R\c C\right)^n$.
The \textbf{ordinary quiver} of $\c C$ has set of vertices $\c C_0$,
and for $x,y\in \c C_0$ the number of arrows $x\to y$ is exactly $dim_k\ _y\c R\c
C_x/_y\c R^2\c C_x$. Finally, we say that $\c C$ is
\textbf{triangular} if its ordinary quiver has no oriented cycle.
All functors are assumed to be $k$-linear functors
between $k$-categories.\\
A functor $F\colon \c E\to \c B$ is called a \textbf{covering functor} (see
\cite{bongartz_gabriel}) if the following properties are satisfied:\\
\indent{a)} $F^{-1}(x)\neq\emptyset$ for any $x\in \c B_0$,\\
\indent{b)}  for any $x_0,y_0\in\c C$  and any $\hat{x}_0,\hat{y}_0\in\c
  E_0$ such that $F(\hat{x}_0)=x_0$
and $F(\hat{y}_0)=y_0$, the following maps induced by
  $F$ are isomorphisms:
$$
\displaystyle{\bigoplus_{F(\hat{y})=y}}\,_{\hat{y}}\c E_{\hat{x}_0}\to\,_{y_0}\c B_{x_0}\ \text{and}\ 
\displaystyle{\bigoplus_{F(\hat{x})=x}}\,_{\hat{y}_0}\c E_{\hat{x}}\to\,_{y_0}\c B_{x_0}\ .
$$
In particular, if $u\in\ _{y_0}\c B_{x_0}$, the inverse images of $u$
by these isomorphisms will be called \textbf{the lifting of $u$}
(w.r.t. $F$) with
source (resp. target) $\hat{x}_0$ (resp. $\hat{y}_0$).
Recall that if $\c E$ is locally
bounded (resp. connected) then so is $\c B$.\\
A \textbf{$G$-category} is a $k$-category $\c C$
endowed with a group morphism $G\to Aut(\c C)$. Moreover, if the
induced action of $G$ on $\c C_0$ is free, then $\c C$ is called a
\textbf{free $G$-category}. The \textbf{quotient category} $\c C/G$ of a free
$G$-category $\c C$ (see \cite{cibils_marcos} for instance) 
has set of objects $\c C_0/G$. For any
$\alpha,\beta\in\c \c C_0/G$ we set:
$$
_{\beta}(\c C/G)_{\alpha}=\left(\bigoplus_{x\in \alpha,y\in\beta}\ _y\c C_x\right)/G 
$$
and the composition is induced by the composition in $\c
C$. The natural projection $\c C\to \c C/G$ is a covering functor.
A \textbf{Galois covering with group $G$} is a functor $F\colon \c
E\to \c B$ with $\c E$ a free
 $G$-category and such that there exists a commutative diagram:\\
\null\hfill\xymatrix{
&\c E \ar@{->}[dl] \ar@{->}[dr]^F& \\
\c E/G \ar@{->}[rr]^{\sim} && \c B
}\hfill\null\\
where $\c E\to \c E/G$ is the natural projection and the horizontal
arrow is an isomorphism. In
particular a Galois covering is a covering functor.
A \textbf{connected Galois covering} is a Galois covering $\c E\to \c
B$ where $\c E$ is connected.\\
A \textbf{$G$-graded category} is a $k$-category $\c C$
such that each morphism space has a decomposition
$
_y\c C_x\ =\oplus_{g\in G}\ _y\c C_x^g
$
satisfying $_z\c C_y^g.\ _y\c C_x^h\subseteq\ _z\c C_x^{gh}$.
The \textbf{smash-product} category (see \cite{cibils_marcos}) $\c C\sharp G$
has set of objects $(\c C\sharp G)_0=\c C_0\times G$, and 
$_{(y,t)}(\c C\sharp G)_{(x,s)}=\ _y\c C_x^{t^{-1}s}$
 for $(x,s)$ and $(y,t)$ in $(\c C\sharp G)_0$. 
The composition in $\c C\sharp G$ is induced by the composition in $\c C$.
The natural projection $F\colon \c C\sharp G\to \c C$, defined by $F(x,s)=x$ and $F(u)=u$ for 
$u\in\ _{(y,t)}(\c C\sharp G)_{(x,s)}\subseteq\ _y\c C_x$, is a Galois
covering with group $G$.
It has been shown in \cite{cibils_marcos} that if  $p\colon \c E\to \c
B$ is  a Galois covering 
with group $G$, then $\c B$ is a $G$-graded category and there
exists a commutative diagram:\\
\null\hfill
\xymatrix{
\c E \ar@{->}[rr]^{\sim}_{\varphi} \ar@{->}[rd]_p & & \c B\sharp G \ar@{->}[ld]\\
&\c B &
}
\hfill\null\\
where $\c B\sharp G\to \c B$ is the natural projection and $\varphi$
is an isomorphism.
\vskip 5pt
\noindent{}\textbf{Quivers with admissible relations}\\
Let $Q$ be a \textbf{locally finite} quiver with set of vertices $Q_0$, set of arrows $Q_1$
and source and target maps $s,t\colon Q_1\to Q_0$ respectively. Recall
that locally finite means that
$s^{-1}(x)$ and $t^{-1}(x)$ are finite sets for any $x\in Q_0$. For
simplicity we will write $x^+$ (resp. $x^-$) for the set $s^{-1}(x)$
(resp. $t^{-1}(x)$).
A (non trivial) \textbf{oriented path} in $Q$ is a non empty sequence
$\alpha_1,\ldots,\alpha_n$ of arrows of $Q$
such that $s(\alpha_{i+1})=t(\alpha_i)$ for any $1\leqslant i\leqslant
n-1$. Such a path is written $\alpha_n\ldots \alpha_1$, its source (resp. target) is
$s(\alpha_1)$ (resp. $t(\alpha_n)$). For each $x\in Q_0$ we will write $e_{x}$
for the (trivial) path of length $0$ and with source and target equal to $x$. 
The \textbf{path category} $kQ$ has set of objects $Q_0$, the
morphism space
$_ykQ_x$ is the vector 
space with basis the set of oriented paths in $Q$ with source $x$ and
target $y$ (including $e_x$ in case $x=y$). The composition
of morphisms in $kQ$ is induced by the concatenation of paths. Notice
that $kQ$ is a free $k$-category in the following sense: for any
$k$-category $\c C$, a functor $kQ\xrightarrow{F}\c C$ is uniquely
determined by the family of
morphisms $\{F(\alpha)\in\ _{F(y)}\c C_{F(x)}\ |\
x\xrightarrow{\alpha}y\in Q_1\}$.
We will denote by $kQ^+$ the ideal of $kQ$ generated by $Q_1$.
Notice also that if $Q_0$ is finite then $kQ$ is also a $k$-algebra,
$kQ=\oplus_{x,y}\ _ykQ_x$, with unit $1=\sum_{x\in Q_0} e_x$, and $kQ^+$
becomes an ideal of this $k$-algebra.
If $r\in\ _ykQ_x$ we define the \textbf{support of $r$} (denoted by
$supp(r)$) to be the set of paths in $Q$ which appear in $r$ with a
non zero coefficient. Moreover, we
define \textbf{normal form of $r$} as an equality of the type $r=\sum_i\lambda_i\
u_i$
such that $\lambda_i\in k^*$ for any $i$ and where the paths $u_i$ are pairwise distinct.
An \textbf{admissible ideal} of $kQ$
is an ideal $I\subseteq kQ$ such that $I\subseteq (kQ^+)^2$ and such that for 
any $x\in Q_0$ there exists $n\geqslant 2$ such that $I$
contains all the paths with length at least $n$ and with source or target $x$.
The couple $(Q,I)$ is then called a \textbf{quiver
with admissible relations} and the quotient category $kQ/I$ is locally bounded.
When $Q_0$ is finite, an admissible ideal $I$ of $kQ$ is exactly an
ideal $I$ of the $k$-algebra $kQ$ such that $(kQ^+)^n\subseteq
I\subseteq (kQ^+)^2$ for some integer $n\geqslant 2$.
Recall from \cite{bongartz_gabriel} that if $\c C$ is a locally bounded $k$-category
then there exists an admissible ideal $I$ for the ordinary quiver $Q$ of
$\c C$ and there exists an isomorphism $kQ/I\xrightarrow{\sim}\c C$.
 Such an isomorphism is called a \textbf{presentation of $\c C$ with
 quiver and (admissible) relations} (or an \textbf{admissible
 presentation} for short). Similarly, if $A$ is a basic finite dimensional
$k$-algebra, an admissible presentation of $A$ is an
isomorphism of $k$-algebras $kQ/I\xrightarrow{\sim}A$ where $(Q,I)$ is
a bound quiver.
\vskip 5pt
\noindent{}\textbf{Transvections, dilatations}\\
A \textbf{bypass} (see \cite{assem_castonguay_marcos_trepode}) in $Q$
is a couple $(\alpha,u)$ where $\alpha\neq u$, $\alpha\in Q_1$ and 
$u$ is a path in $Q$ parallel to $\alpha$ (this means that
$\alpha$ and $u$ share the same source and the same target). A \textbf{double
bypass} is a $4$-tuple $(\alpha,u,\beta,v)$ such that $(\alpha,u)$ and
$(\beta,v)$ are bypasses and such that the arrow $\beta$ appears in
the path $u$. Notice that if $\alpha,\beta$ are distinct parallel arrows of
$Q$, then $(\alpha,\beta,\beta,\alpha)$ is a double bypass. Notice
also that if $u=va$ is an oriented cycle in $Q$ with first arrow $a$,
then $(a, au, a, au)$ is a double bypass. Hence, if $Q$ has no
double bypass, then $Q$ has no distinct parallel arrows and no
oriented cycle. If $A$ is a basic $k$-algebra with quiver $Q$, we
will say that \textbf{$A$ has no double bypass} if $Q$ has
no double bypass. A \textbf{transvection} is an automorphism
of the $k$-category $kQ$ of the form $\varphi_{\alpha,u,\tau}$ where
$(\alpha,u)$ is a bypass, $\tau\in k$ and $\varphi_{\alpha,u,\tau}$ is
given by
$\varphi_{\alpha,u,\tau}(\alpha)=\alpha+\tau\ u$ and
$\varphi_{\alpha,u,\tau}(\beta)=\beta$ for any arrow
$\beta\neq\alpha$ (this uniquely defines $\varphi_{\alpha,u,\tau}$
since $kQ$ is a free $k$-category). Notice that $Q$ has no double bypass if and only
if any two transvections commute. A \textbf{dilatation} is an automorphism $D\colon
kQ\xrightarrow{\sim}kQ$ such that $D(\alpha)\in k^*\alpha$ for any
arrow $\alpha$. Notice that the definitions of transvections and
dilatations are analogous to those of transvection and dilatation
matrices (see \cite[Chap. XIII, § 9]{lang} for instance). Recall that a dilatation
matrix of $GL_n(k)$ is a diagonal invertible matrix with at most one
diagonal entry different from $1$ and a transvection matrix is a matrix with
diagonal entries equal to $1$ and which has at most one non diagonal
entry different from $0$.
\vskip 5pt
\noindent{}\textbf{Fundamental group, coverings of quivers with
relations}\\
Let $(Q,I)$ be a quiver with admissible relations. For each arrow
$\alpha\in Q_1$ we will write $\alpha^{-1}$ for its formal inverse with source
(resp. target) $s(\alpha^{-1})=t(\alpha)$
(resp. $t(\alpha^{-1})=s(\alpha)$). A walk
is an unoriented path in $Q$. More precisely it is a formal product
$u_n\ldots u_1$ of
arrows and of formal inverses of arrows
 such that
$s(u_{i+1})=t(u_i)$ for
any $1\leqslant i\leqslant n-1$. Let $r=t_1u_1+\ldots+t_nu_n \in\ _yI_x$ where 
$t_i\in k^*$ and where the paths $u_i$ are distinct. Then $r$ is
called a \textbf{minimal relation} if $n\geqslant 1$ and if for 
any non empty proper subset $E$ of $\{1,\ldots,n\}$ we have $\sum_{i\in E}t_iu_i\not\in\,_yI_x$.
With this definition, any $r\in I$ can be written as the sum of
minimal relations with pairwise disjoint supports.
Notice that in this definition we do not assume $n$ to be greater than
or equal to $2$ as
done usually (see \cite{martinezvilla_delapena}). This change is done
for simplicity and does not affect the
constructions which follow. 
The \textbf{homotopy relation} of $(Q,I)$
is the smallest equivalence relation $\sim_I$
on the set of walks (in $Q$) which is compatible with the concatenation of walks and such that:\\
\indent{.} $\alpha\alpha^{-1}\sim_I e_y$ and $\alpha^{-1}\alpha\sim_I e_x$ 
for any arrow $x\xrightarrow{\alpha}y$,\\
\indent{.} $u_1\sim_I u_2$ for any minimal relation $t_1u_1+\ldots+t_nu_n$.\\
Notice that in order to compute $\sim_I$ we may restrict ourselves to
any set of minimal relations generating the ideal $I$ (see
\cite{farkas_green_marcos}).
Assume that $Q$ is connected (i.e. $Q$ is connected as an unoriented graph) and let $x_0\in Q_0$.
The \textbf{fundamental group} (see \cite{martinezvilla_delapena}) $\pi_1(Q,I,x_0)$
of $(Q,I)$ at $x_0$ is the set of $\sim_I$-classes of walks starting
and ending at $x_0$. The composition is induced by the 
concatenation of walks and the unit is the $\sim_I$-class of $e_{x_0}$.
Since different choices for $x_0$ give rise to isomorphic fundamental groups (since $Q$ is connected) we shall
write $\pi_1(Q,I)$ for short.
\begin{example}
\label{exple1}(see \cite{assem_delapena})
Assume that $Q$ is the following quiver:
\\
\null\hfill
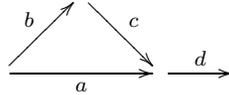
\xymatrix{
 & \ar@{->}[rd]^c & &\\
\ar@{->}[rr]_a \ar@{->}[ru]^b & & \ar@{->}[r]^d & 
}
\hfill\null\\
and set $I=<\ da\ >$ and $J=<\ da-dcb\ >$.
 Then $kQ/I\simeq kQ/J$ whereas $\pi_1(Q,I)\simeq\mathbb{Z}$ and $\pi_1(Q,J)= 0$.
\end{example}
\noindent{A} \textbf{covering} $(Q',I')\xrightarrow{p}(Q,I)$  of quivers
with admissible relations
(see \cite{martinezvilla_delapena}) 
is a quiver morphism $Q'\xrightarrow{p} Q$ such that $p(I')\subseteq I$
and such that:\\
\indent{a)} $p^{-1}(x)\neq\emptyset$ for any $x\in Q_0$,\\
\indent{b)} $x^+\xrightarrow{p} p(x)^+$ and $x^-\xrightarrow{p} p(x)^-$ are bijective for any $x\in Q_0'$,\\
\indent{c)} for any minimal relation $r\in\,_yI_x$ and for any $x'\in p^{-1}(x)$
 there exist $y'\in p^{-1}(y)$ and $r'\in\,_{y'}I'_{x'}$ such that $p(r')=r$,\\
\indent{d)} same statement as c) after interchanging $x$ and $y$.\\
Recall that the automorphism group $Aut(Q,I)$ of a bound quiver $(Q,I)$ is the
group of automorphisms $g\colon Q\xrightarrow{\sim}
Q$ of the quiver $Q$ such that $g(I)\subseteq I$. Assume that
$(Q',I')\xrightarrow{p}(Q,I)$ is a covering, then the \textbf{group of
  automorphisms} of $p$ is defined by $Aut(p)=\{g\in
Aut(Q',I')\ |\ p\circ g=p\}$. If
$(Q',I')\xrightarrow{p}(Q,I)$ is a covering and if $G$ is a subgroup of
$Aut(p)$, then $p$ is called a \textbf{Galois covering with group $G$}
 if $Q$ and $Q'$ are connected and if $G$ acts transitively on
 $p^{-1}(x)$ for any $x\in Q_0$.
Notice that if $(Q',I')\xrightarrow{p}(Q,I)$ is a covering (resp. a Galois covering with
group $G$) then the induced functor
$kQ'/I'\xrightarrow{\bar{p}}kQ/I$ is a covering functor (resp. a Galois
covering with group $G$).
Let $(Q,I)$ be a connected quiver with admissible relations and let $x_0\in Q_0$.
The \textbf{universal cover of $(Q,I)$} is a Galois covering $(\tilde{Q},\tilde{I})\xrightarrow{\pi}(Q,I)$
with group $\pi_1(Q,I,x_0)$ as defined in \cite{martinezvilla_delapena}. One can describe it as follows: $\tilde{Q}_0$ is the set of
$\sim_I$-classes $[w]$ of walks $w$ starting at $x_0$. The arrows of $\tilde{Q}$ are the couples $(\alpha,[w])$
where $\alpha\in Q_1$ and $[w]\in \tilde{Q}_0$ are such that
$s(\alpha)=t(w)$. The source (resp. target) of the arrow $(\alpha,[w])$ is $[w]$ (resp. $[\alpha w]$).
The map $\tilde{Q}\xrightarrow{p}Q$ is defined by $p([w])=t(w)$ and
$p(\alpha,[w])=\alpha$. The ideal $\tilde{I}$ 
is equal to $p^{-1}(I)$. Finally, the action of $\pi_1(Q,I)$ on
$(\tilde{Q},\tilde{I})$ is the following: if $g\in\pi_1(Q,I)$ we may
write $g=[\gamma]$ with $\gamma$ some walk with source and target equal to
$x_0$. Then for any $[w]\in \tilde{Q}_0$ (resp. $(\alpha,[w])\in
\tilde{Q}_1$) we have $g.[w]=[w\gamma^{-1}]$ (resp. $g.(\alpha,[w])=(\alpha,[w\gamma^{-1}])$).
\vskip 5pt
\noindent{}\textbf{Some linear algebra}\\
We introduce here some notions that will be useful in the sequel and
that will be used without reference.  
Let $E$ be a finite dimensional $k$-vector space with
 basis $(e_1,\ldots,e_n)$ and let $(e_1^*,\ldots,e_n^*)$
be the corresponding dual basis (i.e.
$e_i^*(e_i)=1$ and $e_i^*(e_j)=0$ if $j\neq i$).
If $\{r_t\}_{t\in T}$
is a family in $E$, then $Span(r_t\ ;\ t\in T)$ will denote the subspace of $E$ generated by
this family.
 If $r\in E$ we will write $supp(r)$ (the support of
$r$) for the set of
those $e_i's$ appearing in $r$ with a non zero coefficient. Therefore $e_i\in supp(r)$ is
equivalent to $e_i^*(r)\neq 0$. Let $F\subseteq E$ be a
subspace. A non zero element $r\in F$ is called minimal
if it
cannot be written as the sum of two non zero elements of $F$ with
 disjoint supports. We will denote by $\equiv_F$
 the smallest equivalence relation
 on $\{e_1,\ldots,e_n\}$ such that $e_i\equiv_F e_j$ for any
$r\in F$ minimal and any $e_i,e_j\in supp(r)$.
Like in the situation of the homotopy relation of a bound quiver, the
equivalence relation $\equiv_F$ is determined by any generating family
of $F$ made of minimal elements.
Notice that if $E$ is the vector space with basis the set of oriented paths in
 a finite quiver $Q$ (without oriented cycle) and if $I$ is an admissible ideal of $kQ$, then for any
 paths $u$ and $v$ we have: $u\equiv_I v\Rightarrow u\sim_I v$. The
 converse is usually false as one can see in Example~\ref{exple1}
 where $a\sim_J cb$ and $a\not\equiv_J cb$. 
Assume now that the basis of $E$ is totally ordered: $e_1<\ldots<e_n$.
A \textbf{Gröbner} basis of $F$ is a basis $(r_1,\ldots,r_t)$
of $F$ such that:\\
\indent{.} for any $j$ there is some $i_j$ such that $r_j\in
 e_{i_j}+Span(e_i\ ;\ i<i_j)$.\\
\indent{.} $e_{i_j}\not\in supp(r_{j'})$ unless $j=j'$.\\
\indent{.} if $r=e_l+\sum_{i<l}\tau_i\ e_i\in F$ then $e_l=e_{i_j}$
for some $j$.\\
With this definition, $F$ has a unique Gröbner basis which has a
 natural total order: $r_1<\ldots<r_t$ if we assume that
 $i_1<\ldots<i_t$. Moreover, $r_1,\ldots,r_t$ are minimal elements of
 $F$. This last property implies in particular that
$e_i\equiv_F e_j$ if and only if there exists a sequence of integers
$m_1,\ldots,m_p$ such that $e_i\in supp(r_{m_1})$, $e_j\in
supp(r_{m_p})$ and $supp(r_{m_j})\cap supp(r_{m_{j+1}})\neq\emptyset$
for each $j$.
Notice that our definition of Gröbner basis is slightly different from
the classical one (see e.g. \cite{adams_loustaunau}) since we do not
use any multiplicative structure. Moreover, our definition is linked
with the notion of reduced echelon form matrix (see
\cite[p. $65$]{hohn}). We may also point out that
a study of Gröbner bases in
path algebras of quivers has been made in 
\cite{farkas_feustel_green}.\\
We end this paragraph with a reminder on the exponential and on the
logarithm of an endomorphism. If $u\colon E\to E$ is a nilpotent
endomorphism, we define the
\textbf{exponential} of $u$ to be $exp(u)=\sum_{l\geqslant
0}\frac{1}{l!}\ u^l$. Thus, $exp(u)\colon E\to E$ is a
well defined linear isomorphism such that $exp(u)-Id$ is nilpotent.  If $v\colon E\to E$ is an isomorphism such that $v-Id$ is nilpotent, we define the
\textbf{logarithm} of $v$ to be $log(v)=\sum_{l\geqslant
0}(-1)^{l+1}\frac{1}{l}\ (v-Id)^l$. Recall that if $u\colon
E\to E$ is a nilpotent endomorphism, then $log(exp(u))=u$.

\section{Proof of Theorem~\ref{A}}
\label{s1}
In this section we provide the proof of Theorem~\ref{A}
(see also \cite[Thm 1.1]{lemeur}). We fix $A$ a basic connected finite
dimensional $k$-algebra with quiver $Q$. \textbf{Throughout this
  section we will assume that $Q$ has no oriented cycle}.
The proof of Theorem~\ref{A} decomposes into $4$ steps as follows, and we
will devote a subsection to each step:\\
\indent{a)} If $kQ/I$ and $kQ/J$ are isomorphic to $A$
as $k$-algebras, then there exists $\varphi\colon
kQ\xrightarrow{\sim}kQ$ a product of transvections and of a dilatation
such that $\varphi(I)=J$.\\
\indent{b)} If $\varphi(I)=J$ and if $\varphi$ is a
dilatation then $\pi_1(Q,I)\simeq\pi_1(Q,J)$. If $\varphi$ is a
transvection, then there exists a surjective group morphism
$\pi_1(Q,I)\to\pi_1(Q,J)$ or $\pi_1(Q,J)\to\pi_1(Q,I)$, induced by
the identity map on the walks in $Q$.\\
\indent{c)} The homotopy relations $\sim_I$ of the
admissible presentations  $kQ/I$ of $A$ can be displayed as the
vertices of a quiver $\Gamma$ such that for any arrow
$\sim_I\to\sim_J$ the identity map on walks induces a surjective group
morphism $\pi_1(Q,I)\twoheadrightarrow\pi_1(Q,J)$.\\
\indent{d)} If $k$ has characteristic zero and if $Q$ has no double
bypass, then the quiver $\Gamma$ has a unique source. Moreover,
if $\sim_{I_0}$ is the source of $\Gamma$ then $I_0$ fits Theorem~\ref{A}.
\subsection{\normalsize Different presentations of an algebra are linked by products of
  transvections and dilatations}
In order to consider $A$ as a $k$-category we need to choose a
decomposition of the unit into a sum of primitive orthogonal 
idempotents. The following proposition shows that this choice is unrelevant and that we may
fix these idempotents once and for all. We will omit the proof which
is basic linear algebra.
%
%
%
%
\begin{proposition}{\cite[3.1]{lemeur}}
\label{1.1}
Let $I$ and $J$ be admissible
ideals of $kQ$. If $kQ/I\simeq kQ/J$ as $k$-algebras
then there exists $\varphi\colon
kQ\xrightarrow{\sim}kQ$ an automorphism extending
the identity map on $Q_0$ and such that $\varphi(I)=J$.
\end{proposition}
Recall that $GL_n(k)$ is generated by transvection and dilatation matrices.
The following proposition states an analogous result for the group of
automorphisms of $kQ$ extending the identity
map on $Q_0$.
%
%
%
%
\begin{proposition}
\label{1.3}
Let $\c G$ be the
group of automorphisms of $kQ$
 extending the identity map on $Q_0$. Let $\c D\subseteq\c G$ be the
subgroup of the dilatations of $kQ$ and let $\c T\subseteq\c G$ be the
subgroup generated by the transvections. Then $\c T$ is a normal
subgroup and $\c G=\c D\c T=\c T\c D$.
\end{proposition}
\begin{remark} The group of automorphisms of an algebra has already
been studied. More precisely the reader can find in
  \cite{guilasensio_saorin}, \cite{strametz1} and \cite{strametz2} a
  study of the group of outer automorphisms of an algebra.
\end{remark}
\noindent{\textbf{Proof of Proposition~\ref{1.3}:}}
 For any transvection $\varphi=\varphi_{\alpha,u,\tau}$ and any
dilatation $D$ we have $D\varphi
D^{-1}=\varphi_{\alpha,u,\frac{\tau\lambda}{\mu}}$ where $\lambda\in k^*$ and
$\mu\in k^*$ are such that $D(u)=\lambda\ u$ and $D(\alpha)=\mu\
\alpha$. Hence, in order to prove the proposition, it is enough to
prove that $\c G=\c T\c D$. If $\psi\in\c G$ we shall write $n(\psi)$ for the number of
arrows $\alpha\in Q_1$ such that $\psi(\alpha)\not\in k^*\alpha$. Notice
that $n(\psi)=0$ if and only if $\psi\in\c D$. Let us prove by
induction on $n\geqslant 0$ that $R_n:$
``$n(\psi)\leqslant n\Rightarrow\psi\in \c T\c D$'' is true. Obviously $R_0$ is true. Let $n\geqslant
1$, assume that $R_{n-1}$ is true, and let $\psi\in\c G$ such that
$n(\psi)=n$. Hence, there exists $x\xrightarrow{\alpha_1}y\in
Q_1$ such that $\psi(\alpha_1)\not\in k^*\alpha_1$. Let
$\alpha_1,\ldots,\alpha_d$ be the arrows $x\to y$ of $Q$ and let
$E=Span(\alpha_1,\ldots,\alpha_d)\simeq\ _y\left(kQ^+/(kQ^+)^2\right)_x$. Since $kQ\xrightarrow{\psi}kQ$ is
an automorphism,
the composition $f\colon E\hookrightarrow\ _ykQ_x\xrightarrow{\psi}\
_ykQ_x\twoheadrightarrow E$ of $\psi$ with the natural inclusion
 and the natural projection
 is a $k$-linear isomorphism hence an
element of $GL_d(k)$. Thus (see \cite[Chap. XIII Prop. 9.1]{lang}) there exist
transvections matrices $f_1,\ldots,f_l\in GL_d(k)$ such that $f_1\ldots
f_lf(\alpha_i)\in k^*\alpha_i$ for each $i\in\{1,\ldots,d\}$. For each
$f_j$, let $\bar{f}_j\colon kQ\to kQ$ be the automorphism such that
$\bar{f}_j(\alpha_i)=f_j(\alpha_i)$ for each $i\in\{1,\ldots,d\}$ and
such that $\bar{f}_j(\beta)=\beta$ for any arrow $\beta$ not parallel
to $\alpha_1$. In particular, $\bar{f}_j$ is a transvection with
respect to some $\alpha_{i_j}$. Let $g_1=\bar{f}_1\ldots
\bar{f}_l\in\c T$.
Then, $g_1\psi(\alpha_i)\in k^*\alpha_i+(kQ^+)^2$ and if $\beta\in Q_1$
is not parallel to $\alpha_1$
and satisfies $\psi(\beta)\in k^*\beta$ then $g_1\psi(\beta)\in
k^*\beta$. Let $\psi_1=g_1\psi$. By construction, for each
$i\in\{1,\ldots,d\}$, we have
$\psi_1(\alpha_i)=\lambda_i\ \alpha_i+\sum_{j=1}^{n_i}\ \tau_{i,j}\
u_{i,j}$ with $u_{i,j}$ paths of length at least $2$. Let
$\varphi_{i,j}$ be the transvection
$\varphi_{\alpha_i,u_{i,j},-\tau_{i,j}/\lambda_i}$ for each
$i\in\{1,\ldots,d\}$ and each $j\in\{1,\ldots,n_i\}$, and let $g_2\in \c
T$ be
the product of the $\varphi_{i,j}$'s (for any $i\in\{1,\ldots,d\}$ and any $j\in\{1,\ldots,n_i\}$). It is easy to
check that the transvections  $\varphi_{i,j}$ commute between each other so that the
definition of $g_2$ is unambiguous. Since $Q$ has no oriented cycle, we
have $g_2\psi_1(\alpha_i)=\lambda_i\alpha_i$ for each $i$, and $g_2\psi_1(\beta)\in
k^*\beta$ if $\beta\in Q_1$ is not parallel to $\alpha_1$ and satisfies
$\psi_1(\beta)\in k^*\beta$. In particular: $n(g_2g_1\psi)<n(\psi)=n$.
Since $R_{n-1}$ is true, $g_2g_1\psi$ lies in $\c T\c D$ and so does
$\psi$ (recall $g_1,g_2\in\c T$).
 Hence, $R_n$ is true. This achieves the proof of
 Proposition~\ref{1.3}.
\null\hfill$\square$\\
\begin{remark} 
\label{1.2}
Proposition~\ref{1.1} and
Proposition~\ref{1.3} imply that if $I$ and $J$ are admissible ideals of $kQ$
such that $kQ/I\simeq kQ/J$ as $k$-algebras, then there exist
$\varphi_1,\ldots,\varphi_n$ (resp. $\varphi_1',\ldots,\varphi_m'$) a sequence of transvections of $kQ$,
together with $D$ a dilatation such that
$J=D\varphi_n\ldots\varphi_1(I)$ (resp. $J=\varphi_m'\ldots\varphi_1'D(I)$).
\end{remark}
\subsection{\normalsize Comparison of the fundamental group of two presentations
  of an algebra linked by a transvection or a dilatation}
If $I$ is an ideal and $\varphi$ is a dilatation or a transvection,
then $I$ and $\varphi(I)$ are similar enough in order to compare the
associated homotopy relations. Before stating this comparison we prove
two useful lemmas. We fix $I$ an admissible ideal of
$kQ$, we fix $\varphi=\varphi_{\alpha,u,\tau}$ a transvection
($\tau\neq 0$) and we
set $J=\varphi(I)$.
%
%
%
%
\begin{lemma}
\label{1.4}
 Assume that $\alpha\not\sim_I u$
and let $r\in\,_yI_x$ be a minimal relation with normal form 
$r=\sum_C\lambda_c\ \theta_c+\sum_B\lambda_b\ v_b\alpha u_b$ such that
$\alpha$ does not appear in the path $\theta_c$ for any $c\in C$. Then there exists
a minimal relation $r'\in\,_yJ_x$ with normal form $r'=\sum_C\lambda_c\ \theta_c
+\sum_B\lambda_b\ v_b\alpha u_b+\sum_{B'}\lambda_b\tau\ v_buu_b$ where $B'\subseteq B$.
\end{lemma}
\noindent{\textbf{Proof:}} Let us assume that $B\neq\emptyset$ (if $B=\emptyset$, the conclusion is immediate).
Since $Q$ has no oriented cycle, the paths $v_b$ and $u_b$ do not
contain $\alpha$.
Since $r$ is a minimal relation of
 $I$  and since $\alpha\not\sim_I u$, we have $\theta_c\neq v_buu_b$
 for any $c\in C, b\in B$. Therefore, $\varphi(r)$ has
  a normal form $\varphi(r)=\sum_C\ \lambda_c\ \theta_c+
\sum_B\lambda_b\ v_b\alpha u_b+\sum_B\lambda_b\tau\  v_bu u_b\in\,_yJ_x\backslash\{0\}$. Thus
there exists a minimal relation $r'\in\,_yJ_x$ with normal form
$r'=\sum_{C'}\lambda_c\ \theta_c+
\sum_{B_1'}\lambda_b\  v_b\alpha u_b+\sum_{B'}\lambda_b\tau\  v_bu u_b$
such that $\emptyset\neq B_1'\subseteq B$, $C'\subseteq C$ and
 $B'\subseteq B$.
Hence $\varphi^{-1}(r')$ has a normal form $\varphi^{-1}(r')=\sum_{C'}\lambda_c\ \theta_c+
\sum_{B_1'}\lambda_b\  v_b\alpha u_b+\sum_{B'\backslash B_1'}\lambda_b\tau\  v_bu u_b-
\sum_{B_1'\backslash B'}\lambda_b\tau\  v_bu u_b\in\,_yI_x\backslash\{0\}$. Since
$r\in\,_yI_x$ is a minimal relation and since $\alpha\not\sim_I u$ we infer that
there exists a minimal relation $r''\in\,_yI_x$ with normal form
$r''=\sum_{C''}\lambda_c\ \theta_c+\sum_{B''}\lambda_b\  v_b\alpha u_b$
such that $C''\subseteq C'\subseteq C$ and $\emptyset\neq B''\subseteq B_1'$. This forces $C''=C$
and $B''=B$ because $r\in\,_yI_x$ is a minimal relation. Thus $C'=C$ and
$B_1'=B$. Hence we have a minimal relation $r'\in\,_yJ_x$ with normal form
$r'=\sum_C\lambda_c\ \theta_c
+\sum_B\lambda_b \ v_b\alpha u_b+\sum_{B'}\lambda_b\tau\  v_buu_b$ as announced.\hfill$\square$
%
%
%
%
\begin{lemma}
\label{1.5}
 Assume that $\alpha\sim_Ju$ and let $r\in\,_yI_x$ be a 
minimal relation. Then $v\sim_Jw$
for any $v,w\in supp(r)$.
\end{lemma}
\noindent{\textbf{Proof:}} We may write $r=\sum_C\lambda_c\ \theta_c+\sum_B\lambda_b\  v_b\alpha u_b
+\mu_b\ v_buu_b$ where:\\
. $\lambda_c,\lambda_b\in k^*$ and $\mu_b\in k$ 
for any $c\in C$ and $b\in B$,\\
. the paths $\theta_c$, $v_b\alpha u_b$, $v_{b'}uu_{b'}$ ($c\in C,b,b'\in B$)
are pairwise distinct,\\
. for any $c\in C$, the path $\theta_c$ does 
not contain $\alpha$.\\
Hence $\varphi(r)=\sum_C\lambda_c\ \theta_c+\sum_B\lambda_b\  v_b\alpha u_b+(\mu_b+\tau\lambda_b)\ v_buu_b\in\,_yJ_x$
and there exists a decomposition $\varphi(r)=r_1+\ldots +r_n$ 
where $r_i\in\,_yJ_x$ is a minimal relation
and $supp(r_i)\cap supp(r_j)=\emptyset$ if $i\neq j$. 
If $B=\emptyset$ then $\varphi(r)=r\in\,_yJ_x$
is a minimal relation and the lemma is proved. Hence we
 may assume that $B\neq\emptyset$. This implies that
for any $i\in\{1,\ldots,n\}$ there exists $b\in B$ such that
 $v_b\alpha u_b\in supp(r_i)$ or $v_buu_b\in supp(r_i)$
(if this is not the case then $r_i=\sum_{C'}\lambda_c\ \theta_c$ 
for some non empty subset $C'$ of $C$, thus
$\varphi^{-1}(r_i)=\sum_{C'}\lambda_c\,\theta_c\in\,_yI_x$ 
which contradicts the minimality of $r$).
Let $\equiv$ be the smallest equivalence relation on the 
set $\{1,\ldots,n\}$ such that:
$i\equiv j$ if there exists $b\in B$ such that
$v_b\alpha u_b\in supp(r_i)$ and $v_buu_b\in supp(r_j)$.
Since the $r_i$'s are minimal relations of $J$ and since $\alpha\sim_Ju$,
 we get: if $i\equiv j$ then $v\sim_J w$ 
for any $v,w\in supp(r_i)\sqcup supp(r_j)$. 
Let $\c O\subseteq\{1,\ldots,n\}$ be a $\equiv$-orbit
 and let $r'=\sum_{i\in \c O}r_i\in\,_yJ_x$.
Hence  $r'=\sum_{C'}\lambda_c\ \theta_c+\sum_{B'}\lambda_b\  v_b\alpha
u_b+(\mu_b+\tau\lambda_b)\ v_buu_b$
where $C'\subseteq C$ and
$\emptyset\neq B'\subseteq B$.
This implies that 
$\varphi^{-1}(r')=\sum_{C'}\lambda_c\ \theta_c+\sum_{B'}\lambda_b\  v_b\alpha u_b+\mu_b\ v_buu_b\in\,_yI_x$
and the minimality of $r$ yields $C'=C$, $B'=B$, $r'=\varphi(r)$ 
and $\c O=\{1,\ldots,n\}$. Hence $\{1,\ldots,n\}$ is an
$\equiv$-orbit. Therefore
$v\sim_Jw$ for any $v,w\in supp(\varphi(r))$. And since $\alpha\sim_Ju$ we infer that $v\sim_J w$
 for any $v,w\in supp(r)$.\hfill$\square$
%
%
%
%
\vskip 5pt
We can state the announced comparison now. For short,
the word \textit{generated} stands for:
\textit{generated as an equivalence relation compatible with the
concatenation of walks and such that $\alpha^{-1}\alpha\sim_I e_x$,
$\alpha\alpha^{-1}\sim_I e_y$ for any arrow $x\xrightarrow{\alpha}y$}.
\begin{proposition}{\cite[3.2]{lemeur}}
\label{1.6}
Let $I$ be an
admissible ideal of $kQ$, let $\varphi$ be
an automorphism of $kQ$ and set $J=\varphi(I)$. If $\varphi$ is a dilatation,
then $\sim_I$ and $\sim_J$ coincide. Assume now that
$\varphi=\varphi_{\alpha,u,\tau}$ is a transvection.\\
\indent{a)} if $\alpha\sim_Iu$ and $\alpha\sim_Ju$ then $\sim_I$ and
$\sim_J$ coincide.\\
\indent{b)} if $\alpha\not\sim_Iu$ and $\alpha\sim_J u$ then $\sim_J$
is generated by $\sim_I$ and  $\alpha\sim_Ju$.\\
\indent{c)} if $\alpha\not\sim_Iu$ and $\alpha\not\sim_Ju$ then $I=J$
and $\sim_I$ and $\sim_J$ coincide.
\end{proposition}
\begin{remark}  The
  following implication (symmetric to $b)$):\\
\null\hfill
if $\alpha\sim_Iu$ and $\alpha\not\sim_Ju$ then $\sim_I$ is generated
  by $\sim_J$ and $\alpha\sim_Iu$\hfill\null\\
is also satisfied since
  $\varphi_{\alpha,u,\tau}^{-1}=\varphi_{\alpha,u,-\tau}$
\end{remark}
\noindent{\textbf{Proof of Proposition~\ref{1.6}:}} If $\varphi$ is a dilatation, then $\sim_I$
and $\sim_J$ coincide because for any $r\in\ _ykQ_x$ we have
$supp(r)=supp(\varphi(r))$ and because $r$ is a minimal relation of $I$ if and
only if the same holds for $\varphi(r)$ in $J$. Let us assume that
$\varphi=\varphi_{\alpha,u,\tau}$ is a transvection.\\
a) Lemma~\ref{1.5} applied to
$I,J,\varphi$ (resp. $J,I,\varphi^{-1}=\varphi_{\alpha,u,-\tau}$) shows that any two paths appearing in a same minimal
relation of $I$  (resp. $J$) are $\sim_J$-equivalent
(resp. $\sim_I$-equivalent). Hence $\sim_I$ and $\sim_J$ coincide.\\
b) Let $\equiv$ be the equivalence relation generated by:
$(v\sim_Iw\Rightarrow v\equiv w)$ and $\alpha\equiv u$. Our aim is to
show that $\sim_J$ and $\equiv$ coincide. Thanks to Lemma~\ref{1.5} we
have: $v\equiv w\Rightarrow v\sim_J w$.
Let $Min(I)$ be the set of the minimal relations of $I$. For each $r\in
Min(I)$ let us fix a normal form $r=\sum_C\lambda_c\
\theta_c+\sum_B\lambda_b\ v_b\alpha u_b$
satisfying the hypotheses of Lemma~\ref{1.4}. Hence there exists
$B'\subseteq B$ and a minimal
relation $r_1$ of $J$ with normal form $r_1=\sum_C\lambda_c\ \theta_c
+\sum_B\lambda_b\ v_b\alpha u_b+\sum_{B'}\lambda_b\tau\ v_buu_b$.
Thus $\varphi(r)-r_1=\sum_{B\backslash
B'}\tau\lambda_b\ v_buu_b\in J$ can be written as a sum
$r_2+\ldots+r_{n_r}$ of minimal relations of $J$
with pairwise disjoint supports. In particular,
$\varphi(r)=r_1+\ldots+r_{n_r}$ where each $r_i\in J$ is a minimal
relation. Notice that any two paths appearing in
$r_1$ are $\equiv$-equivalent because of the normal form of $r_1$ and
because of the definition of
$\equiv$. With these notations, the set $\{r_i\
|\ r\in Min(I)\ and\ 1\leqslant i\leqslant n_r\}$ is made of minimal
relations of $J$ and generates the ideal $J$. Thus, in order to show
that $\sim_J$ and $\equiv$ coincide, it is enough to show that any two
paths appearing in some $r_i$ are $\equiv$-equivalent.
Let $r\in Min(I)$, let $i\in\{1,\ldots,n_r\}$, and let $v,w\in
supp(r_i)$. We have already proved that if $i=1$ then
$v\equiv w$, thus we may assume that $i\geqslant 2$. Keeping the above
notations for the normal form of $r$,  there exist $b,b'\in B$
such that $v=v_buu_b$ and $w=v_{b'}uu_{b'}$. Since $\alpha\equiv u$
and since any two paths appearing in $r_1$ are $\equiv$-equivalent we get
$v=v_buu_b\equiv v_b\alpha u_b\equiv v_{b'}\alpha u_{b'}\equiv v_{b'}uu_{b'}=w$.
Hence any two paths appearing in some $r_i$ are $\equiv$-equivalent.
This implies that $\sim_J$ and $\equiv$ coincide. Therefore, $\sim_J$ is generated by $\sim_I$ and
$\alpha\sim_J u$.\\
c) Let $r\in I$ be a minimal relation of $I$ and apply
Lemma~\ref{1.4} to $r$. Since $\alpha\not\sim_J u$, we infer that
$r\in 
J$. Since $I$ is generated by its minimal relations we get $I\subseteq
J$. Finally, $I=J$ because $I$ and $J$ have the same dimension.\hfill$\square$
\begin{remark}
\label{1.12}
In the situation b) of Proposition~\ref{1.6}, the identity map on the walks
of $Q$ induces a surjective group morphism $\pi_1(Q,I)\twoheadrightarrow\pi_1(Q,J)$.
\end{remark}
Proposition~\ref{1.6} allows us to prove the following
result which has already been proved in \cite{bardzell_marcos}. Recall that the algebra $kQ/I$,
where $I$ is admissible, is
called \textbf{constricted} if $dim\ _y(kQ/I)_x=1$ for any arrow $x\to y$ of
$Q$.
\begin{proposition}[see also \cite{bardzell_marcos}]
\label{bardzell_marcos}
Assume that $A$ is constricted. Then different admissible
presentations of $A$ yield the same
homotopy relation.
In particular, they have isomorphic fundamental groups.
\end{proposition}
\noindent{\textbf{Proof:}}  Let $kQ/I\simeq A$ be any admissible presentation.
If $(\alpha,u)$ is a bypass in $Q$ then $u\in I$
because $A$ is constricted and $I$ is admissible.
In particular, for any $\tau\in k$ and any $r\in kQ$ we have: $\varphi_{\alpha,u,\tau}(r)-r$ belongs to the ideal
 of $kQ$ generated by $u\in I$ and therefore
 $\varphi_{\alpha,u,\tau}(r)-r\in I$. This shows that
 $\varphi(I)\subseteq I$ and that $\varphi(I)=I$ ($I$ is finite
 dimensional because $Q$ has no oriented cycle)
 for any transvection $\varphi$. Let $kQ/J\simeq A$ be another
 admissible presentation.
From Remark~\ref{1.2} we know that there exist a dilatation $D$ and
transvections $\varphi_1,\ldots,\varphi_n$ such that
$J=D\varphi_n\ldots\varphi_1(I)$. We deduce from what we have proved
above that $\varphi_n\ldots\varphi_1(I)=I$. Hence,
$J=D(I)$. Proposition~\ref{1.6} implies that $\sim_I$ and $\sim_J$
coincide. Therefore, $\pi_1(Q,I)$ and $\pi_1(Q,J)$ are isomorphic.
\hfill$\square$
\vskip 5pt
If $\sim$ and $\sim'$ are homotopy relations, we will say that $\sim'$
is a \textbf{direct successor} (see also \cite[Sect. 3]{lemeur}) of $\sim$ if there exist admissible
ideals $I$ and $J$ of $kQ$, together with a transvection
$\varphi=\varphi_{\alpha,u,\tau}$ such that $\sim=\sim_I$,
$\sim'=\sim_J$, $J=\varphi(I)$, $\alpha\not\sim_Iu$ and
$\alpha\sim_Ju$. Notice that $I,J,\varphi$ need not be unique.
\subsection{\normalsize The quiver $\Gamma$ of the homotopy relations of the
  presentations of the algebra}
\begin{definition}{\cite[4.1]{lemeur}}
  We define the quiver $\Gamma$ as
  follows:\\
  . $\Gamma_0$ is the set of homotopy relations of the
  admissible presentations of $A$:
$$
\Gamma_0=\{\ \sim_I\ |\ I\ \text{is admissible and}\ kQ/I\simeq A\}
$$
  . there is an arrow $\sim\to\sim'$ if and only if $\sim'$ is
  a direct successor of $\sim$.
\end{definition}
\begin{example}
  Assume that $A=kQ/I$ where $Q$ is
\\
\null\hfill
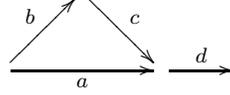
\xymatrix{
 & \ar@{->}[rd]^c & &\\
\ar@{->}[rr]_a \ar@{->}[ru]^b & & \ar@{->}[r]^d & 
}
\hfill\null\\
  and $I=<da>$. Let $J=<da-dcb>$. Using Proposition~\ref{1.6} one can
  show that $\Gamma$ is equal to:
  \xymatrix{\sim_I\ar@{->}[r] & \sim_J}.
  Notice that the identity map on walks induces
  a surjective group morphism
  $\mathbb{Z}\simeq\pi_1(Q,I)\twoheadrightarrow\pi_1(Q,J)\simeq 1$.
\end{example}
The author thanks Mariano Su\'arez-Alvarez for the following remark:
\begin{remark}
  A homotopy relation is determined by its restriction to the paths in
  $Q$ with length at most the radical length of $A$. Thus there are only
  finitely many homotopy relations. This argument shows that $\Gamma$ is finite.
\end{remark}
%
%
%
%
The following proposition states some additional properties of $\Gamma$ and is a
direct consequence of Remark~\ref{1.2} and Proposition~\ref{1.6}.
\begin{proposition}
\label{1.7}
Assume that $Q$ has no oriented cycle and let $m$ be
the number of bypasses in $Q$. Then $\Gamma$ is connected and has no
oriented cycle. Any vertex of $\Gamma$ is the source of at most $m$
arrows and any oriented path in $\Gamma$ has length at most $m$.
\end{proposition}
\begin{remark}
\label{1.13}
 According to Remark~\ref{1.12}, if there is a path in
  $\Gamma$ with source $\sim_I$ and target $\sim_J$, then the identity
  map on the walks in $Q$ induces a surjective group morphism
  $\pi_1(Q,I)\twoheadrightarrow\pi_1(Q,J)$.
Moreover, since $\Gamma$ is finite, any vertex of $\Gamma$ is the target of a
(finite) path whose source is a source of $\Gamma$  (i.e. a
vertex with no arrow ending at it). As a consequence, if $\Gamma$
has a unique source $\sim_{I_0}$, then the fundamental group of any admissible
presentation of $A$ is a quotient of $\pi_1(Q,I_0)$.
\end{remark}
\subsection{\normalsize The uniqueness of the source of  $\Gamma$ and the proof  of
  Theorem~\ref{A}}
Notice that until now we have used neither the characteristic of 
$k$ nor the possible non existence of a double bypass in $Q$.
These hypotheses will be needed in order to prove the uniqueness of
the source of $\Gamma$. The complete proof of the uniqueness of the
source of $\Gamma$
is somewhat technical. For this reason we deal with the technical
considerations in the two lemmas that follow.
%
%
%
%
\begin{lemma}
\label{1.8}
Let $E$ be a finite dimensional $k$-vector space endowed with a
totally ordered basis $e_1<\ldots<e_n$. Assume that $k$ has
characteristic zero. Let $\nu\colon E\to E$ be
a linear map such that $\nu(e_i)\in Span(e_j\ ;\ j<i)$ for any
$i\in\{1,\ldots,n\}$, and let $I$ and $J$ be two subspaces of $E$ such
that the following conditions are satisfied:\\
\indent{a)} $\psi(I)=J$ where $\psi\colon E\to E$ is equal to
$exp(\nu)$.\\
\indent{b)} if $e_i\in supp(\nu(e_j))$ then
$e_i\not\equiv_I e_j$ and $e_i\not\equiv_J e_j$.\\
Then $I$ and $J$ have the same Gröbner basis and $I=J$.
\end{lemma}
\noindent{\textbf{Proof:}} Let us prove Lemma~\ref{1.8} by induction
on $n$. If $n=1$ the equality is obvious so let us assume that $n>1$
and that the conclusion of Lemma~\ref{1.8} holds for dimensions less
than $n$. We will denote by $r_1<\ldots<r_p$
(resp. $r_1'<\ldots<r_p'$) the Gröbner basis of $I$ (resp. of $J$) and
we will write $i_1,\ldots,i_p$ (resp. $i_1',\ldots,i_p'$) for the
integers such that $r_j\in e_{i_j}+Span(e_i\ ;\ i<i_j)$ (resp. $r_j'\in e_{i_j'}+Span(e_i\ ;\ i<i_j')$).
In order to prove that  $I=J$ we will prove the four following facts:
\begin{enumerate}
\item[a)] the two sequences $i_1<\ldots<i_p$ and
  $i_1'<\ldots<i_p'$ coincide,
\item[b)] $\psi(r_1)=r_1'$,
\item[c)] $r_1=r_1'$ and $\nu(r_1)=0$ (using the induction hypothesis on $E/k.e_1$),
\item[d)] $r_2=r_2',\ldots,r_p=r_p'$ (using the induction hypothesis on $E/k.r_1$).
\end{enumerate}
\vskip 5pt
\noindent{a)} For simplicity let us set $E_i=Span(e_j\ ;\ j\leqslant i)$. Since
$\nu(e_j)\in E_{j-1}$ and $r_j\in e_{i_j}+E_{i_j-1}$, and since
$\psi=exp(\nu)$, we get $\psi(r_j)\in
J\cap\left(e_{i_j}+E_{i_j-1}\right)$ for any $j$. Hence, the
definition of the Gröbner basis of $J$ forces 
$\{i_1,\ldots,i_p\}\subseteq\{i_1',\ldots,i_p'\}$ and the cardinality
and the ordering on these two sets imply that $i_1=i_1',\ldots,i_p=i_p'$
\vskip 5pt
\noindent{b)} Since $i_1=i_1'$ we infer that 
$\psi(r_1)-r_1'\in J\cap E_{i_1-1}$. Then, the definition of the
Gröbner basis of $J$ forces $\psi(r_1)-r_1'=0$.
\vskip 5pt
\noindent{c)} Let us prove that $r_1=r_1'$. 
Notice that the definition of a Gröbner basis and the equalities
$\psi(r_1)=r_1'$ and $\psi(e_1)=e_1$ force:
$e_1\in
I\Leftrightarrow r_1=e_1\Leftrightarrow r_1'=e_1\Leftrightarrow e_1\in
J$. Hence we may assume that
$e_1\not\in I$ and $e_1\not\in J$.\\
Let $\widetilde{E}=E/k.e_1$ and let $\pi\colon E\twoheadrightarrow
\widetilde{E}$ be the natural projection. We will write $\widetilde{x}$
for $\pi(x)$. Similarly we set $\widetilde{I}=\pi(I)$ and $\widetilde{J}=\pi(J)$. In
particular $\widetilde{E}$ has a totally ordered basis:
$\widetilde{e}_2<\ldots<\widetilde{e}_n$. Since $\nu(e_1)=0$ and since
$\psi(e_1)=e_1$, the mappings $\nu$ and $\psi$ induce linear mappings
$\widetilde{\nu},\widetilde{\psi}\colon\widetilde{E}\to\widetilde{E}$. It follows from the
properties of $\nu$ and $\psi$ that $\widetilde{\psi}(\widetilde{I})=\widetilde{J}$, that $\widetilde{\nu}(\widetilde{e}_i)\in Span(\widetilde{e}_j\ ;\ 2\leqslant
j<i)$ for any $i\geqslant 2$, that
$\widetilde{\psi}=exp(\widetilde{\nu})$,
and that $supp(\widetilde{\nu}(\widetilde{e}_i))=\{\ \widetilde{e}_j|\ j\geqslant 2\ and\
e_j\in supp(\nu(e_i))\}$ for any $i\geqslant 2$.
Moreover, with the definition of the Gröbner basis of $I$ we get:
\begin{enumerate}
\item[.] $\widetilde{r}_j\in \widetilde{e}_{i_j}+Span(\widetilde{e}_i\ ;\
i<i_j)$ for any $j$ (recall that $e_1\not\in I$),
\item[.] $supp(\widetilde{r}_j)=\{\widetilde{e}_i\ |\ i\geqslant 2\
and\ e_i\in supp(r_j)\}$ for any $j$.
\end{enumerate}
 Therefore $\widetilde{r}_1<\ldots<\widetilde{r}_p$
is the Gröbner basis of $\widetilde{I}$ and: 
$\widetilde{e}_i\equiv_{\widetilde{I}}\widetilde{e}_j\Rightarrow e_i\equiv_I e_j$.
Similarly $\widetilde{r}_1'<\ldots<\widetilde{r}_p'$
is the Gröbner basis of $\widetilde{J}$ and: $\widetilde{e}_i\equiv_{\widetilde{J}}\widetilde{e}_j\Rightarrow e_i\equiv_J e_j$.
Using the above description of $supp(\widetilde{\nu}(\widetilde{e}_i))$ together with
the above link between $\equiv_I$ (resp. $\equiv_J$) and
$\equiv_{\widetilde{I}}$ (resp. $\equiv_{\widetilde{J}}$) we infer that:\\
$$
\widetilde{e}_i\not\equiv_{\widetilde{I}}\widetilde{e}_j\ \text{and}\ \widetilde{e}_i\not\equiv_{\widetilde{J}}\widetilde{e}_j\ \text{as
soon as}\ \widetilde{e}_j\in supp(\widetilde{\nu}(\widetilde{e}_i))
$$
For this reason we may apply the induction hypothesis to $\widetilde{E}$,
$\widetilde{I}$ and $\widetilde{J}$. Hence $\widetilde{I}$ and $\widetilde{J}$ have
the same Gröbner basis and $\widetilde{r}_1=\widetilde{r}_1'$ i.e. $r_1'=r_1+\lambda\ e_1$ with $\lambda\in k$.
 Therefore $(\psi-Id)(r_1)=\lambda\ e_1$, and since $\psi(e_1)=e_1$ we
 get $\nu(r_1)=log(\psi)(r_1)=\lambda\ e_1$.
Assume that $\lambda\neq 0$ i.e. $e_1\in supp(\nu(r_1))$. Thus there exists $e_i\in supp(r_1)$
such that $e_1\in supp(\nu(e_i))$. This implies that $e_1\not\equiv_I
e_i$, and since any two elements in
$supp(r_1)$ are $\equiv_I$-equivalent, this forces $e_1\not\in
supp(r_1)$.
Hence $e_i,e_1\in supp(r_1')=supp(r_1)\sqcup\{e_1\}$ and therefore
$e_i\equiv_J e_1$. This contradicts $e_1\in supp(\nu(e_i))$ and shows
that $\lambda=0$, that $r_1=r_1'$ and that $\nu(r_1)=0$.
\vskip 5pt
\noindent{d)} Let us show that $r_2=r_2',\ldots,r_p=r_p'$.
For this purpose we will apply the induction hypothesis to $\bar{E}=E/k.r_1$. Let $q\colon
E\twoheadrightarrow \bar{E}$ be the natural projection. 
We will write $\bar{e}_i$ (resp. $\bar{I}$, $\bar{J}$, $\bar{r}_j$,
$\bar{r}_j'$) for $q(e_i)$ (resp. $q(I)$, $q(J)$, $q(r_j)$, $q(r_j')$).
Hence $\bar{E}$ has a totally ordered basis:
$\bar{e}_1<\ldots<\bar{e}_{i_1-1}<\bar{e}_{i_1+1}<\ldots<\bar{e}_n$.
Since $\nu(r_1)=0$ and since $\psi(r_1)=r_1$, the mappings $\nu$ and
$\psi$ induce linear mappings $\bar{\nu},\bar{\psi}\colon \bar{E}\to\bar{E}$. These mappings obviously satisfy
$\bar{\psi}(\bar{I})=\bar{J}$, $\bar{\nu}(\bar{e}_i)\in Span(\bar{e}_j\ ;\ j\neq i_1\ and\
j<i)$ for any $i\neq i_1$, and $\bar{\psi}=exp(\bar{\nu})$.
Moreover, our choice for the basis of $\bar{E}$ and the definition of the Gröbner basis of $I$
imply that:
\begin{enumerate}
\item[.] $supp(\bar{r}_j)=\{\bar{e}_i\ |\ e_i\in supp(r_j)\}$ for any
$j\geqslant 2$,
\item[.] $\bar{r}_2<\ldots<\bar{r}_p$ is the Gröbner basis of
$\bar{I}$.
\end{enumerate}
These two properties imply in particular that:
$\bar{e}_i\equiv_{\bar{I}}\bar{e}_j\Rightarrow e_i\equiv_Ie_j$ for any
$i,j\neq i_1$. The corresponding properties hold for $\bar{J}$
(replace $r_j$ by $r_j'$, $I$ by $J$ and $\bar{I}$ by
$\bar{J}$). Thus, in order to apply the induction
hypothesis to $\bar{E}$ it only remains to prove that: $\bar{e}_j\in
supp(\bar{\nu}(\bar{e}_i))\Rightarrow
\bar{e}_i\not\equiv_{\bar{I}}\bar{e}_j\ \text{and}\
\bar{e}_i\not\equiv_{\bar{J}}\bar{e}_j$ for any $i,j\neq i_1$. Assume
that $i,j\neq i_1$ satisfy $\bar{e}_j\in
supp(\bar{\nu}(\bar{e}_i))$. From the definition of $\bar{E}$ and
$\bar{\nu}$ we know that:
\begin{enumerate}
\item[.] $supp(\bar{\nu}(\bar{e}_i))=\{\bar{e}_l\ |\ e_l\in supp(\nu(e_i))\}$ if $e_{i_1}\not\in
supp(\nu(e_i))$,
\item[.] $supp(\bar{\nu}(\bar{e}_i))\subseteq \{\bar{e}_l\ |\ e_l\in
supp(\nu(e_i))\ \text{and}\ l\neq i_1\}\cup \{\bar{e}_l\ |\ l<i_1\ \text{and}\ e_l\in
supp(r_1)\}$ if $e_{i_1}\in supp(\nu(e_i))$.
\end{enumerate}
Let us distinguish the cases $e_j\in supp(\nu(e_i))$ and $e_j\not\in supp(\nu(e_i))$:\\
\indent{$\cdot$} if $e_j\in supp(\nu(e_i))$ then $e_i\not\equiv_I e_j$ and
$e_i\not\equiv_Je_j$ and the above comparison between $\equiv_I$
(resp. $\equiv_J$) and $\equiv_{\bar{I}}$
(resp. $\equiv_{\bar{J}}$) yields
$\bar{e}_i\not\equiv_{\bar{I}}\bar{e}_j$ and
$\bar{e}_i\not\equiv_{\bar{J}}\bar{e}_j$.\\
\indent{$\cdot$} if $e_j\not\in supp(\nu(e_i))$ then necessarily $e_{i_1}\in
supp(\nu(e_i))$ and $e_j\in supp(r_1)$. Since $r_1=r_1'$, the property
 $e_j\in supp(r_1)$ implies that  $e_j\equiv_I
e_{i_1}$ and $e_j\equiv_Je_{i_1}$. On the other hand, the property $e_{i_1}\in
supp(\nu(e_i))$ implies that 
 $e_{i_1}\not\equiv_I e_i$ and
$e_{i_1}\not\equiv_J e_i$. Therefore $e_j\not\equiv_Ie_i$ and
$e_j\not\equiv_Je_i$ and finally $\bar{e}_j\not\equiv_{\bar{I}}\bar{e}_i$ and
$\bar{e}_j\not\equiv_{\bar{J}}\bar{e}_i$.\\
Thus all the conditions of Lemma~\ref{1.8} are satisfied for
$\bar{E},\bar{I},\bar{J},\bar{\nu}$. For this reason we can
apply the induction hypothesis which gives: $\bar{I}$ and $\bar{J}$
have the same Gröbner basis. We infer that
$q(r_i)=q(r_i')$ for each $i=2,\ldots,p$. Hence for each $i\geqslant
2$ there exists $\lambda_i\in k$ such that $r_i=r_i'+\lambda_i\ r_1$, and
$\lambda_i$ is necessarily zero because
$e_{i_1}^*(r_i)=e_{i_1}^*(r_i')=0$ (cf the definition of a Gröbner basis).
Therefore $r_i=r_i'$ for each $i=1,\ldots,p$ and $I=J$
as announced.\hfill$\square$\\
%
%
%
%
\begin{lemma}
\label{1.9}
Let
$\varphi\colon kQ\to kQ$ be an automorphism extending the identity map
on $Q_0$. Let $I$ be an admissible ideal of $kQ$ and set $J=\varphi(I)$.
Suppose that $k$ has characteristic zero. Suppose that for any arrow $\alpha$
there is a normal form $\varphi(\alpha)=\alpha+\sum_i\lambda_i\ u_i$
where each $u_i$ satisfies: $\alpha\not\sim_I u_i$,
 $\alpha\not\sim_J u_i$ and $\varphi(a)=a$ for any arrow appearing
in $u_i$ (in particular $\varphi(u_i)=u_i$). Then $I$ and $J$ coincide.
\end{lemma}
\noindent{\textbf{Proof:}} Let $E$ be the vector space $kQ=\oplus_{x, y}\
_ykQ_x$. Hence $E$ is finite dimensional since $Q$ has no oriented
cycle, and $I$ and $J$ can be considered as subspaces of $E$. In order
to apply Lemma~\ref{1.8} to $E,I,J$, we need to exhibit a totally
ordered basis of $E$ together with a mapping $\nu\colon E\to E$. Let us
take the family of paths in $Q$ for the basis of $E$. The following
construction of a total order $<$ on this basis is taken from
\cite{farkas_feustel_green}. Let us fix a total order on $Q_0\cup Q_1$ (which
is finite) and let $\prec$ be the induced lexicographical order on the
paths in $Q$ ($e_x\prec u$ if $u$ is non trivial). If $u$ is a path we let $W(u)$ be the number of arrows
$\alpha\in Q_1$ appearing in $u$ and such that
$\varphi(\alpha)\neq\alpha$. Hence, for any $\alpha\in Q_1$, we have
$W(\alpha)=0$ if $\varphi(\alpha)=\alpha$ and $W(\alpha)=1$ if
$\varphi(\alpha)\neq\alpha$. The total order $<$ is then defined as
follows:
$$
u<v\ \Leftrightarrow\left\{
\begin{array}{lll}
W(u)<W(v)\\
or\\
W(u)=W(v)\ \text{and}\ u\prec v
\end{array}\right.
$$
This yields: $e_1<\ldots< e_n$ a totally
ordered basis of $E$ made of the paths in $Q$. Notice that
with this basis, the equivalence relations $\equiv_I$ and $\sim_I$
(resp. $\equiv_J$ and $\sim_J$) satisfy the following property:
$e_i\equiv_Ie_j\Rightarrow e_i\sim_I e_j$ (resp. $e_i\equiv_Je_j\Rightarrow
e_i\sim_J e_j$). Let $\nu\colon kQ\to kQ$ be the derivation (i.e. the
$k$-linear map such that $\nu(vu)=\nu(v)u+v\nu(u)$ for any $u$ and $v$) such
that $\nu(e_x)=0$ for any $x\in Q_0$ and $\nu(\alpha)=\varphi(\alpha)-\alpha$ for any arrow $\alpha\in
Q_1$. Thus, for any path $u$ and any $v\in supp(\nu(u))$ there exist an
arrow $\alpha\in Q_1$ together with paths $u_1,u_2,u_3$ such that
$u=u_3\alpha u_1$, $v=u_3u_2u_1$ and $u_2\in
supp(\nu(\alpha))$. Notice that with the assumptions made on $\varphi$,
this implies that $e_i\not\equiv_I e_j$ and $e_i\not\equiv_J e_j$ as
soon as $e_j\in supp(\nu(e_i))$. Moreover, for any $\alpha\in Q_1$ and
any $u\in supp(\nu(\alpha))$ we have $W(u)=0$, hence $\nu\circ
\nu(\alpha)=0$. Since $\nu\colon kQ\to kQ$ is a derivation, we infer
that: $e_j\in supp(\nu(e_i))\Rightarrow W(e_j)<W(e_i)\Rightarrow
e_j<e_i$. Hence $\nu(e_i)\in Span(e_j\ ;\ j<i)$ for any $i$. In order
to apply Lemma~\ref{1.8}, it only remains to prove that
$J=exp(\nu)(I)$. To do this it suffices to prove that
$\varphi=exp(\nu)$. 
Since $\nu$ is a derivation, $exp(\nu)\colon kQ\to kQ$ is an
automorphism such that $exp(\nu)(e_x)=e_x$ for any $x\in Q_0$ (recall that $\nu(e_x)=0$). Moreover, if $\alpha\in Q_1$ then $\nu^2(\alpha)=0$ and
$\nu(\alpha)=\varphi(\alpha)-\alpha$, therefore
$exp(\nu)(\alpha)=\varphi(\alpha)$. Hence $\varphi$ and $exp(\nu)$ are
automorphisms of $kQ$ and they coincide on $Q_0\cup Q_1$. This implies that $\varphi=exp(\nu)$.
Hence, the data
$E,I,J,\nu$ together with the ordered basis $e_1<\ldots<e_n$ satisfy
the hypotheses of Lemma~\ref{1.8} which implies that
$I=J$.\hfill$\square$
\vskip 5pt
The uniqueness of the source of $\Gamma$ is given by the following result.
%
%
%
%
\begin{proposition}{\cite[4.3]{lemeur}}
\label{1.10}
Assume that $A$ satisfies the hypotheses made before stating Theorem~\ref{A}, then $\Gamma$ has a unique source.
\end{proposition}
\noindent{\textbf{Proof:}} Notice that any two transvections of $kQ$
 commute since $Q$ has no double bypass. Let $\sim$ and $\sim'$ be
 sources of $\Gamma$. Let $I$ and $J$ be admissible ideals of $kQ$
 such that $kQ/I\simeq A\simeq kQ/J$ and such that $\sim=\sim_I$ and
 $\sim'=\sim_J$. According to Remark~\ref{1.2} there exist a sequence of transvections
 $\varphi_1=\varphi_{\alpha_1,u_1,\tau_1},\ldots,\varphi_n=\varphi_{\alpha_n,u_n,\tau_n}$
 of $kQ$ and a dilatation $D$ such that
 $J=\varphi_n\ldots\varphi_1D(I)$. Thanks to Proposition~\ref{1.6} we know
 that $\sim_I=\sim_{D(I)}$. Thus, in order to prove
 that $\sim=\sim'$, we may assume that $D=Id_{kQ}$ and
 $J=\varphi_n\ldots\varphi_1(I)$. Moreover we may assume that $n$ is
 the smallest non negative integer such that there exist $I$, $J$ and
 a sequence of transvections $\varphi_1,\ldots,\varphi_n$ satisfying
 $\sim=\sim_I$, $\sim'=\sim_J$ and
 $J=\varphi_n\ldots\varphi_1(I)$. Let us prove that
 $\alpha_i\not\sim_I u_i$ for any $i\in\{1,\ldots,n\}$. If $i$ is such
 that $\alpha_i\sim_I u_i$ then Proposition~\ref{1.6} implies that
 $\sim_I=\sim_{\varphi_i(I)}$ since $\sim_I$ is a source of
 $\Gamma$. Hence $\sim=\sim_{\varphi_i(I)}$, $\sim'=\sim_J$ and
 $J=\varphi_n\ldots\varphi_{i+1}\varphi_{i-1}\ldots\varphi_1(\varphi_i(I))$
 which contradicts the minimality of $n$. Thus $\alpha_i\not\sim_I
 u_i$ for any $i$ and the same arguments apply to $J$ since
 $I=\varphi_1^{-1}\ldots\varphi_n^{-1}(J)$ and since $\sim_J$ is a source of
 $\Gamma$. Hence $\alpha_i\not\sim_J u_i$ for any $i$. This shows that
 the data $I,J,\varphi_n\ldots\varphi_1$ satisfy the hypotheses of
 Lemma~\ref{1.9}. We infer that $I=J$ and that $\sim=\sim'$. This
 shows that $\Gamma$ has a unique
 source.\hfill$\square$
\vskip 5pt
Proposition~\ref{1.10} and Remark~\ref{1.13} prove
Theorem~\ref{A}:
\begin{THMM}{(see also \cite[Thm 1.1]{lemeur})}
\label{1.11}
Let $A$ be a basic connected finite dimensional algebra over a field $k$
of characteristic zero. If the quiver $Q$ of $A$ has no double
bypass, then there exists a presentation $kQ/I_0\simeq A$ with
quiver and admissible relations such that for any other admissible presentation $kQ/I\simeq
A$, the identity map on walks induces a surjective group morphism $\pi_1(Q,I_0)\twoheadrightarrow \pi_1(Q,I)$.
\end{THMM}
The following example shows that one cannot remove the hypothesis on
the characteristic of $k$ in Proposition~\ref{1.10}:
\begin{example}
Let $Q$ be the following quiver without double bypass:
\\
\null\hfill
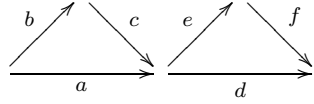
\xymatrix{
 & \ar@{->}[rd]^c & & \ar@{->}[rd]^f & \\
\ar@{->}[ru]^b \ar@{->}[rr]_a & & \ar@{->}[ru]^e \ar@{->}[rr]_d & &
}
\hfill\null\\
Set $u=cb$ and $v=fe$. Set
$A=kQ/I_0$ where $I_0=<d a+vu,va+d u>$. Then
$\pi_1(Q,I_1)=\mathbb{Z}/2$. Let $I_1$ and $I_2$ be the ideals defined
below:\\
\indent{$\bullet$} $I_1=\varphi_{a,u,1}(I_0)=<d a+d u+vu,
va+d u+vu>$,\\
\indent{$\bullet$}
$I_2=\varphi_{a,u,-1}\circ\varphi_{d,v,-1}(I_0)=<d a,va+d
u-2vu>$.\\
Hence $A\simeq kQ/I_1\simeq kQ/I_2$. If $car(k)=0$,  then $\pi_1(Q,I_1)=\pi_1(Q,I_2)=1$ and $\Gamma$ is equal to \xymatrix{\sim_{I_0}\ar@{->}[r]&\sim_{I_1}}.
Suppose now that $car(k)=2$. Then $I_2=<d a,va+d u>$,
$\pi_1(Q,I_0)\simeq\mathbb{Z}/2$, $\pi_1(Q,I_1)=1$, $\pi_1(Q,I_2)\simeq\mathbb{Z}$ and
$\Gamma$ is equal to $\sim_{I_0}\rightarrow\sim_{I_1}\leftarrow\sim_{I_2}$.
Hence $\Gamma$ has two sources. Notice that the identity map on walks
induces a surjective group morphism $\pi_1(Q,I_2)\twoheadrightarrow\pi_1(Q,I_0)$.
Notice also that one can build similar examples for any non zero value $p$ of $car(k)$ by taking for $Q$ a sequence of $p$ bypasses.
\end{example}
\section{Preliminaries on covering functors}
\label{s2}
In this section we give some useful facts on covering
functors.
\begin{lemma}
\label{2.0}
Let $p\colon\c E\to \c B$ and $q\colon\c E'\to \c B$ be covering functors where $\c E$ is connected. Let $r,r'\colon\c E\to
\c E'$ be such that $q\circ r=q\circ r'=p$. If there exists
$x_0\in\c E_0$ such that $r(x_0)=r'(x_0)$ then $r=r'$.
\end{lemma}
\noindent{\textbf{Proof:}} Since $q$ is a covering functor, for any $u\in\ _y\c E_x\backslash\{0\}$ we
have:\\
\null\hfill$(r(x)=r'(x)\ \text{or}\ r(y)=r'(y))\Rightarrow\left(r(u)=r'(u),\ r(x)=r'(x)\ \text{and}\  r(y)=r'(y)\right)$\hfill$(\star)$\\ 
Assume that there exists $x_0\in\c E_0$ such that $r(x_0)=r'(x_0)$.
Since $\c E$ is connected, for any $x\in \c E_0$ there exists a sequence $x_0,\ldots,x_n=x$
of objects of $\c E$ together with a non zero morphism between $x_i$
and $x_{i+1}$ for any $i$. This implies (thanks to $(\star)$) that $r(x)=r'(x)$. Thus
$r$ and $r'$ coincide on $\c E_0$, and $(\star)$ implies $r=r'$.\hfill$\square$
\vskip 5pt
The following proposition generalises the result
\cite[Prop. 3.3]{martinezvilla_delapena}. Using Lemma~\ref{2.0} its
proof is immediate.
\begin{proposition}
\label{2.1}
Let $F\colon \c E\to \c B$ be a covering functor where $\c E$ is
 connected.
 Then $\c E$ is an $Aut(F)$-category. Moreover, $F$ is a Galois
covering if and only if $Aut(F)$ acts transitively on each $F^{-1}(x)$.
Finally, if $F$ is Galois covering with group $G$, then $G=Aut(F)$.
\end{proposition}
\begin{proposition}
\label{2.2}
Let $p\colon \c E\to\c B$ and $q\colon \c F\to \c E$ be functors where
$\c E$ is connected and set $r=p\circ q\colon \c F\to \c B$.
Then $p,q,r$
are covering functors as soon as two of them are so.
\end{proposition}
\noindent{\textbf{Proof:}} Using basic linear algebra arguments, it is
easy to verify the two following facts: if $p$ and $q$ (resp. $q$ and $r$) are covering functors
then so is $r$ (resp. $p$). If $p$ and $r$ are
covering functors, then $q$ satisfies the condition b) in the
definition of a covering functor. We only need to prove that
 $q^{-1}(x)\neq\emptyset$ for any $x\in \c
E_0$. The condition b) implies that $q^{-1}(x)\neq\emptyset\Leftrightarrow
q^{-1}(y)\neq\emptyset$ as soon as $_y\mathcal{E}_x\neq 0$.
Since $\mathcal{E}$ is connected and $q^{-1}(q(x))\neq\emptyset$ for
any $x\in\mathcal{F}_0$ we deduce that $q^{-1}(x)\neq\emptyset$ for any
  $x\in\c E_0$.\hfill$\square$\\
\begin{proposition}
\label{2.3}
Let $p\colon \c C\to \c B$ (resp. $q\colon \c C'\to \c B$) be a
connected Galois
covering with group $G$ (resp. $G'$)
and assume there exists a commutative diagram of $k$-categories and
$k$-linear functors where
$\varphi$ is an isomorphism extending the identity map on $\c B_0$:\\
\null\hfill
\xymatrix{
\c C \ar@{->}[r]^r \ar@{->}[d]_p & \c C' \ar@{->}[d]^q \\
\c B \ar@{->}[r]^{\sim}_{\varphi} & \c B 
}
\hfill\null\\
Then there exists a unique
mapping $\lambda\colon G\to G'$ such that $r\circ g=\lambda(g)\circ r$
for any $g\in G$. Moreover $\lambda$ is a surjective
morphism of groups and $r$ is a Galois covering with group $Ker(\lambda)$.
\end{proposition}
\noindent{\textbf{Proof:}} 
Thanks to Proposition~\ref{2.2}, $r$ is a covering functor.
Fix $\hat{x}_0\in\c C$ and set
$x_0=p(\hat{x}_0)$. For any $g\in Aut(p)$ we have
$q(r(\hat{x}_0))=x_0=q(r(g(\hat{x}_0)))$. Since $q$ is Galois with
group $G'$, there exists a unique $\lambda(g)\in G'$ such that
$\lambda(g)(r(\hat{x}_0))=r(g(\hat{x}_0))$, and Lemma~\ref{2.0} yields
$\lambda(g)\circ r=r\circ g$. Hence: $
\left(\forall g\in G\right)\ \left(\exists !\lambda(g)\in G'\right)\
\lambda(g)\circ r=r\circ g$.
This last property shows the existence and the uniqueness of
$\lambda$. It also shows that $\lambda\colon G\to G'$ is a group
morphism and that $Aut(r)=Ker(\lambda)$. Moreover, $\lambda$ is surjective because of its definition
and because $p$ is Galois with group $G$. Finally,
Proposition~\ref{2.1} shows that $r$ is a Galois covering with group
$Ker(\lambda)$.\hfill$\square$
\section{The universal cover of an algebra}
\label{s3}
In this section we will prove Theorem~\ref{B}. Let $Q$ be
a connected quiver without oriented cycle and fix
 $x_0\in Q_0$ for the computation of the groups $\pi_1(Q,I)$.
If there is no ambiguity we shall write $[w]$ for the homotopy class of
a walk $w$.
%
%
%
%
\begin{lemma}
\label{3.1}
Let $I$ be an admissible ideal of $kQ$, let $D$ be
a dilatation of $kQ$ and set $J=D(I)$. 
Let $\lambda\colon
\pi_1(Q,I)\xrightarrow{\sim}\pi_1(Q,J)$ be the isomorphism given by
Proposition~\ref{1.6}.
 Let $p\colon
(\tilde{Q},\tilde{I})\to (Q,I)$ (resp. $q\colon
(\hat{Q},\hat{J})\to(Q,J)$) be the universal Galois covering with group
$\pi_1(Q,I)$ (resp. $\pi_1(Q,J)$). Then there exists an isomorphism
$\psi\colon
k\tilde{Q}/\tilde{I}\xrightarrow{\sim}k\hat{Q}/\hat{J}$ such that
the following diagram commutes:\\
\null\hfill
\xymatrix{
k\tilde{Q}/\tilde{I} \ar@{->}[r]^{\psi} \ar@{->}[d]_{\bar{p}} &
k\hat{Q}/\hat{J} \ar@{->}[d]^{\bar{q}}\\
kQ/I \ar@{->}[r]^{\bar{D}} & kQ/J
}\hfill\null\\
where $\bar{D},\bar{p}$ and $\bar{q}$ are induced by $D,p$ and $q$ respectively.\\
Moreover, $\psi$ satisfies: $\psi\circ g=\lambda(g)\circ\psi$ for any $g\in \pi_1(Q,I)$.
\end{lemma}
\noindent{\textbf{Proof:}} We have $\hat{Q}=\tilde{Q}$ since $\sim_I$
and $\sim_J$ coincide (see Proposition~\ref{1.6}). Set $\hat{D}\colon
k\tilde{Q}\to k\hat{Q}$ to be defined by: $\hat{D}(a,[w])=(D(a),[w])$
for any arrow $(a,[w])\in\tilde{Q}_1$. By construction $\hat{D}$ is an
automorphism of $kQ$ and
$\hat{D}(\tilde{I})=\hat{J}$. Set $\psi\colon k\tilde{Q}/\tilde{I}\xrightarrow{\sim} 
k\hat{Q}/\hat{J}$ to be induced by $\hat{D}$. Then it is easy to
check all the
announced properties.\hfill$\square$

%
%
%
%
\begin{lemma}
\label {3.2}
Let $I$ be an admissible ideal of $kQ$, let $\varphi=\varphi_{\alpha,u,\tau}$ be
a transvection, set $J=\varphi(I)$ and assume that
$\alpha\sim_Ju$. Let $\lambda\colon
\pi_1(Q,I)\twoheadrightarrow\pi_1(Q,J)$ be the surjection given by
Proposition~\ref{1.6}.
Denote by $p\colon
(\tilde{Q},\tilde{I})\to (Q,I)$ (resp. by $q\colon
(\hat{Q},\hat{J})\to(Q,J)$) the universal Galois covering with group
$\pi_1(Q,I)$ (resp. $\pi_1(Q,J)$). Then there exists a Galois covering
$\psi\colon
k\tilde{Q}/\tilde{I}\xrightarrow{\sim}k\hat{Q}/\hat{J}$ with group
$Ker(\lambda)$ and such that
the following diagram commutes:\\
\null\hfill
\xymatrix{
k\tilde{Q}/\tilde{I} \ar@{->}[r]^{\psi} \ar@{->}[d]_{\bar{p}} &
k\hat{Q}/\hat{J} \ar@{->}[d]^{\bar{q}}\\
kQ/I \ar@{->}[r]^{\bar{\varphi}} & kQ/J
}\hfill\null\\
where $\bar{\varphi},\bar{p}$ and $\bar{q}$ are induced by $\varphi,p$
and $q$ respectively.\\
Moreover, $\psi$ satisfies: $\psi\circ g=\lambda(g)\circ\psi$ for any $g\in \pi_1(Q,I)$.
\end{lemma}
\noindent{\textbf{Proof:}} Let $\varphi'\colon k\tilde{Q}\to k\hat{Q}$
be defined by: $\varphi'([w])=[w]$ for any $[w]\in\tilde{Q}_0$,
$\varphi'(\beta,[w])=(\beta,[w])$ for any $(\beta,[w])\in
\tilde{Q}_1$ such that $\beta\neq\alpha$, and
$\varphi'(\alpha,[w])=(\alpha,[w])+\tau(u,[w])$ for any
$(\alpha,[w])\in\tilde{Q}_1$. Then $\varphi'$ is well defined
since $\alpha\sim_J u$. Moreover, $\varphi\circ
p(a)=q\circ\varphi'(a)$ for any $a\in\tilde{Q}_1$, and
$\varphi'(\tilde{I})\subseteq \hat{J}$. Let $\psi\colon
k\tilde{Q}/\tilde{I}\to k\hat{Q}/\hat{J}$ be induced by $\varphi'$.
Thus $\bar{q}\circ\psi=\bar{\varphi}\circ \bar{p}$. Let
$g=[\gamma]\in\pi_1(Q,I)$ and let $[w]\in\tilde{Q}_0$.
Then $\psi\circ
g([w])=\psi([w\gamma^{-1}])=[w\gamma^{-1}]=\lambda(g)([w])=\lambda(g)\circ
\psi([w])$. The Lemma~\ref{2.0} implies that $\psi\circ
g=\lambda(g)\circ \psi$ for any $g\in\pi_1(Q,I)$. Finally,
Proposition~\ref{2.3} gives: $\psi$ is a Galois covering with group $Ker(\lambda)$.\hfill$\square$
%
%
%
%
\begin{lemma}
\label{3.3}
Let $A$ be a basic and connected finite dimensional $k$-algebra with ordinary quiver
$Q$. Assume that $k$ has characteristic zero and that $Q$ has no
double bypass. Let $\sim_{I_0}$ be the unique source of $\Gamma$ and
 $\sim_I$ be a vertex of $\Gamma$. Then there exist a
sequence
$\varphi_1,\ldots,\varphi_n$ ($\varphi_i=\varphi_{\alpha_i,u_i,\tau_i}$)
of transvections and a dilatation $D$ such that:\\
\indent{a)} $I=D\varphi_n\ldots\varphi_1(I_0)$,\\
\indent{b)} if $I_i$ is the ideal $\varphi_i\ldots\varphi_1(I_0)$ then
$\alpha_i\sim_{I_i}u_i$.
\end{lemma}
\noindent{\textbf{Proof:}} We shall write $[n]$ for the
set $\{1,\ldots,n\}$. Remark~\ref{1.2} implies that
$I=D\psi_1\ldots\psi_m(I_0)$ where
the $\psi_i$'s are transvections and $D$ is a dilatation.
Set $J=D^{-1}(I)=\psi_1\ldots\psi_m(I_0)$. Thus we only need to prove
that the conclusion of Lemma~\ref{3.3} holds for $J$. Let $R_m$ be the
property: \textit{``If $J$ is the image of $I_0$ by a product of $m$
transvections, then there exists a sequence
$\varphi_1,\ldots,\varphi_n$ of transvections such that
$J=\varphi_n\ldots\varphi_1(I_0)$ and which satisfies the property b)
of Lemma~\ref{3.3}''}. Let us prove that $R_m$ is true by induction on $m\geqslant
0$. Obviously $R_0$ is true, let $m\geqslant 1$ and let us assume that
$R_{m-1}$ is true. Let $J=\psi_1\ldots\psi_m(I_0)$ where $\psi_i=\varphi_{a_i,v_i,t_i}$.
 Assume first that there
exists $i_0\in [m]$ such that $a_{i_0}\sim_Jv_{i_0}$.
Set
$J'=\psi_{1}\ldots\psi_{i_0-1}\psi_{i_0+1}\ldots\psi_m(I_0)$. Thanks
to $R_{m-1}$, there exists a sequence $\varphi_1,\ldots,\varphi_n$
of transvections such that $J'=\varphi_n\ldots\varphi_1(I_0)$ and
which satisfies the property b) of Lemma~\ref{3.3}. The sequence
$\varphi_1,\ldots,\varphi_n,\psi_{i_0}$ shows that $R_m$ is true when
such an $i_0$ exists. 
Assume now that for any $i\in[m]$ we have $a_i\not\sim_J
v_i$. Let $\varphi=\psi_m\ldots\psi_1$. Lemma~\ref{1.9}, applied to
the data $I_0,J,\varphi$, shows that $J=I_0$.
Hence $R_m$ is true (with $n=0$) in this situation as well. This
achieves the proof of Lemma~\ref{3.3}.\hfill$\square$
\vskip 5pt
The following proposition shows how a Galois covering of
locally bounded $k$-category is induced by a covering of quivers with
relations. Notice that this proposition makes no assumption on the
ordinary quiver of the involded $k$-categories (in particular, the
quiver may have loops, oriented cycles, multiple arrows...). It
generalises \cite[prop 3.4,
3.5]{martinezvilla_delapena}. The proof uses the
ideas presented in \cite[sect. 3]{green}.
%
%
%
\begin{proposition}
\label{3.4}
Let $F\colon \hat{\c C}\to \c C$ be a Galois covering with group $G$ where
$\c C$ and $\hat{\c C}$ are locally bounded. 
 Then, there
exist admissible presentations $\varphi\colon kQ/I'\xrightarrow{\sim} \c C$
and $\psi\colon k\hat{Q}/\hat{I}\xrightarrow{\sim} \hat{\c C}$ 
and a covering of quiver with relations 
$p\colon(\hat{Q},\hat{I})\to(Q,I')$, such that 
$\varphi$ restricts to the identity map $Q_0=\c C_0\to \c C_0$ on $\c
C_0$ and such that the following diagram is commutative:\\
\null\hfill\xymatrix{
k\hat{Q}/\hat{I} \ar@{->}[r]^{\psi} \ar@{->}[d]_{\bar{p}} & 
\hat{\c C} \ar@{->}[d]^F \\
kQ/I' \ar@{->}[r]^{\varphi}  & \c C
}\hfill\null\\
where $\bar{p}$ is induced by $p$. If $\hat{\c C}$ is
connected, then $p$ is Galois with group $G$.
\end{proposition}
\noindent{\textbf{Proof:}} Using \cite[Thm 3.8]{cibils_marcos} we may assume
that $\c C$ is $G$-graded, that $\c C'=\c C\sharp G$ and that $\c
C'=\c C \sharp G\xrightarrow{F}\c C$ is the natural projection. Since
$\hat{\c C}$ and $\c C$ are locally bounded,
\cite[3.3]{bongartz_gabriel} implies that any morphism in $\c R\c C$
is the sum of images (under $F$) of morphisms in $\c R\hat{\c
  C}$. Since the image under $F$ of a morphism in $\hat{\c C}$ is
a homogeneous morphism in $\c C$, we deduce that 
the ideals $\c R\c
C$ and $\c R^2\c C$ are homogeneous. Thus, for any $x\neq y\in Q_0$ there
exist homogeneous elements  $_yu_x^{(1)},\ldots,\ _yu_x^{(\,_yn_x)}$ of $\ _{y}\c
R\c C_{x}$ giving rise to a basis of
$_{y}\left(\c R\c C/\c R^2\c C\right)_{x}$. In
particular, $_yn_x$ is equal to the number of arrows $x\to
y$ in $Q$. Let $\mu\colon kQ\to \c C$ be defined
as follows: $\mu(x)=x$ for any $x\in Q_0=\c C_0$, and $\mu$ induces
a bijection between the set of arrows $x\to y$ of $Q$ and
$\{_yu_x^{(1)},\ldots,\ _yu_x^{(\,_yn_x)}\}$ for any $x\neq y\in
Q_0$. 
Set $I'=Ker(\mu)$. Hence $I'$ is admissible and $\mu$ induces
an isomorphism $\varphi\colon kQ/I'\xrightarrow{\sim}\c C$.
 The following construction of $p$ uses the ideas of 
Green in \cite[Sect. 3]{green}.
The $k$-category $kQ$ is a $G$-graded as follows: a path $u$ in $Q$ is
homogeneous of degree the degree of $\mu(u)$. By using the $G$-grading on $\c C$, it is easy to
check that $I'$ is homogeneous and that $\varphi\colon kQ/I'\to \c C$
is homogeneous of degree $1_{G}$.
 Let $\hat{Q}$ be the
quiver as follows: $\hat{Q}_0=Q_0\times G$, and the arrows
$(x,s)\xrightarrow{\alpha}(y,t)$ in $\hat{Q}_1$ are exactly the arrows
$x\xrightarrow{\alpha}y$ in $Q_1$ with degree $t^{-1}s$. Let
$p\colon \hat{Q}\to Q$ be defined by: $p(x,s)=x$ and
$p((x,s)\xrightarrow{\alpha}(y,t))=\alpha$ for any $(x,s)\in\hat{Q}_0$ and
any $(x,s)\xrightarrow{\alpha}(y,t)\in\hat{Q}_1$. Let
$\hat{I}\subseteq \hat{Q}$ be the admissible ideal $p^{-1}(I')$ of $k\hat{Q}$. According to
\cite[Sect. 3]{green}, $p$ is a
covering, and if $\hat{Q}$ is connected then $p$ is Galois with group
$G$. In particular $\bar{p}\colon k\hat{Q}/\hat{I}\to kQ/I'$ is a
covering functor. Let $\nu\colon k\hat{Q}\to\c C'=\c C\sharp
G$ be as follows: $\nu(x,s)=(\varphi(x),s)$ for any
$(x,s)\in\hat{Q}_0$, and if
$(x,s)\xrightarrow{\alpha}(y,t)\in\hat{Q}_1$ then $\nu(\alpha)=\mu(p(\alpha))\in\
_{\varphi(y)}\c C_{\varphi(x)}^{t^{-1}s}=\ _{(\varphi(y),t)}\c
C'_{(\varphi(x),s)}$. Therefore $F\circ \nu=\varphi\circ p$, and since $\varphi$ is an
isomorphism, we have $\hat{I}=Ker(\nu)$. Let
$\psi\colon k\hat{Q}/\hat{I}\to \c C'$ be induced by $\nu$. Hence
$\psi\colon \hat{Q}_0\to\hat{\c C}_0$ is bijective,
$\psi$ is faithful and $\varphi\circ
\bar{p}=F\circ\psi$. Moreover $\psi$ is full because $\bar{p}$ and $F$
are covering functors. Thus, $\psi$ is an isomorphism. Finally, if
$\c C'$ is connected then $\hat{Q}$ is
connected and this implies that $p$ is a Galois covering with
group $G$.\hfill$\square$
\vskip 5pt
Recall from \cite{assem_skowronski} that a triangular algebra is called
simply connected if the fundamental group of any
presentation of this algebra is trivial. A triangular
algebra is simply connected if and only if it has no proper connected
Galois covering, see \cite{skowronski}. The following corollary generalises this
characterisation to non triangular algebras, it is a direct
consequence of Proposition~\ref{3.4}
\begin{corollary}
Let $A$ be a basic connected finite dimensional $k$-algebra. Then
the following assertions are equivalent:
\begin{enumerate}
\item[.] for any admissible presentation $kQ/I\simeq A$ we have $\pi_1(Q,I)=1$,
\item[.] $A$ has no proper Galois covering $\c C\to A$ 
 with $\c C$ connected and locally bounded.
\end{enumerate}
\end{corollary}
\begin{remark} Proposition~\ref{3.4} does not necessarily hold
when $F$ is a covering functor and not a Galois covering. As an
example, set $\c C=kQ$ where $Q$ is equal to:\\
\null\hfill
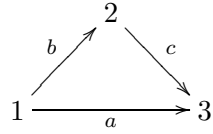
\xymatrix{
& 2\ar@{->}[rd]^c&\\
1\ar@{->}[ru]^b \ar@{->}[rr]_a&&3
}
\hfill\null\\
set $G=\mathbb{Z}/2=<\sigma|\sigma^2>$ and set $\c C'=kQ'$ where $Q'$ is the quiver:\\
\null\hfill
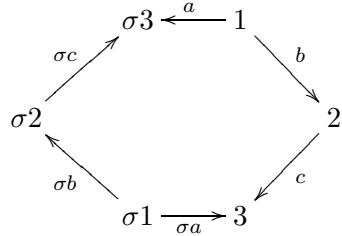
\xymatrix{
& \sigma 3 & 1 \ar@{->}[l]_a \ar@{->}[rd]^b & \\
\sigma 2 \ar@{->}[ru]^{\sigma c} & & & 2\ar@{->}[ld]^c\\
& \sigma 1\ar@{->}[lu]^{\sigma b} \ar@{->}[r]_{\sigma a} & 3 
}\hfill\null\\
Set $F\colon \c C'\to\c C$ to be defined by: $F(b)=F(\sigma b)=b$,
$F(c)=F(\sigma c)=c$, $F(a)=a$ and $F(\sigma a)=a+cb$. Then $F$ is a
covering functor. The group $Aut(F)$ is trivial therefore $F$ is not Galois,
and $F$ cannot be induced by any covering of bound quivers. Notice that if $F\colon \c
C'\to\c C$ is a covering functor and if the ordinary quiver of $\c C$
has no bypass, then $F$
is induced by a covering of bound quivers.
\end{remark}
%
%
%
%
\begin{THMM}
\label{3.5}
Assume that $A$ satisfies the hypotheses made before stating Theorem~\ref{A}.
Let $\varphi_0\colon kQ/I_0\simeq A$ be an admissible presentation
such that $\sim_{I_0}$ is the source of $\Gamma$. Let
$(\tilde{Q},\tilde{I}_0)\xrightarrow{p_0}(Q,I_0)$
be the universal Galois covering with group $\pi_1(Q,I_0)$
and let $k\tilde{Q}/\tilde{I}_0\xrightarrow{\bar{p}_0} kQ/I$ be
induced by $p_0$.
For any connected Galois covering $F\colon \c C'\to A$ with group $G$ there
exist an isomorphism $kQ/I_0\xrightarrow{\sim}A$ equal to $\varphi_0$
on objects, a Galois
covering $F'\colon k\tilde{Q}/\tilde{I}_0\to \c C'$ with group $N$ a
normal subgroup of $\pi_1(Q,I_0)$ such that the following diagram commutes:\\
\null\hfill\xymatrix{
k\tilde{Q}/\tilde{I}_0 \ar@{->}[r]^{F'} \ar@{->}[d]_{\bar{p}_0} &  \c C' \ar@{->}[d]^F \\
kQ/I_0 \ar@{->}[r]^{\sim} & A
}\hfill\null\\
Moreover, there is an exact sequence of groups: $1\to N\to \pi_1(Q,I_0)\to G\to 1$.
\end{THMM}
\noindent{\textbf{Proof:}} Let $\c C'\xrightarrow{F}A$ be a
connected Galois covering with group $G$. Thanks to Proposition
\ref{3.4} we may assume that there exists a Galois covering (with group $G$) of bound quivers
$(Q',I')\xrightarrow{q}(Q,I)$ such that: $A=kQ/I$, 
$\c C'=kQ'/I'$ and $F\colon\c C'\to A$ is induced by $q$. Let
$(\hat{Q},\hat{I})\xrightarrow{p}(Q,I)$ be the universal Galois
covering with group $\pi_1(Q,I)$. Thus (see
\cite{martinezvilla_delapena}) there exists a Galois covering
$(\hat{Q},\hat{I})\xrightarrow{r}(Q',I')$
such that $q\circ r=p$. Hence we have a commutative diagram
(denoted by $\c D$):\\
\null\hfill\xymatrix{
k\hat{Q}/\hat{I} \ar@{->}[d]_{\bar{p}} \ar@{->}[r]^{\bar{r}}
& kQ'/I' \ar@{->}[d]_{\bar{q}} \ar@{->}[r]_{\psi}^{\sim} & \c C'
\ar@{->}[d]^F\\
kQ/I \ar@{=}[r]^{Id}&kQ/I\ar@{->}[r]^{\sim}_{\varphi}&A
}\hfill\null\\
Since $\sim_{I_0}$ is the source of $\Gamma$, Lemma~\ref{3.3}
implies that
there exist both a sequence of transvections
$\varphi_1=\varphi_{\alpha_1,u_1,\tau_1},\ldots,\varphi_n=\varphi_{\alpha_b,u_n,\tau_n}$
of $kQ$ and a dilatation $D$ such that
$I=D\varphi_n\ldots\varphi_1(I_0)$ and such that
$\alpha_i\sim_{I_i}u_i$ if $I_i=\varphi_i\ldots\varphi_1(I_0)$ for any
$i$. Lemma~\ref{3.1} and Lemma~\ref{3.2} applied to $D,I,I_n$
and $\varphi_i,I_{i-1},I_i$ respectively yield the following
commutative diagrams denoted by $\c D'$ and $\c T_i$ respectively:\\
\null\hfill\xymatrix{
kQ^{(n)}/I^{(n)} \ar@{->}[d]^{\bar{p_n}} \ar@{->}[r]&
k\hat{Q}/\hat{I} \ar@{->}[d]^{\bar{p}}&
kQ^{(i-1)}/I^{(i-1)} \ar@{->}[d]^{\bar{p}_{i-1}}\ar@{->}[r]&
kQ^{(i)}/I^{(i)} \ar@{->}[d]^{\bar{p}_i}\\
kQ/I_n \ar@{->}[r]^{\bar{D}}&
kQ/I&
kQ/I_{i-1} \ar@{->}[r]^{\bar{\varphi_i}}&
kQ/I_i&
}\hfill\null\\
where $\bar{\varphi_i}$ (resp. $\bar{D}$) is induced by $\varphi_i$
(resp. $D$) and $kQ^{(i)}/I^{(i)}\xrightarrow{\bar{p}_i}kQ/I_i$ is
induced by the universal Galois covering
$(Q^{(i)},I^{(i)})\xrightarrow{p_i}(Q,I_i)$ with group $\pi_1(Q,I_i)$.
If we connect $\c T_1,\ldots,\c
T_n,\c D'$ and $\c D$ we get the announced commutative diagram.
Finally the announced properties of $F'$ are
given by Proposition~\ref{2.3}.\hfill$\square$
\begin{remark}
Using the universal property in Theorem~\ref{3.5} it is easily verified
 that if there exists a Galois covering $\c C'\to A$ such
that $\c C'$ is simply connected (i.e.  the fundamental group of any
presentation of $\c C'$ is trivial), then $\c C'\simeq
k\tilde{Q}/\tilde{I}_0$.
\end{remark}
 One may wish to use the more general framework of Galois
categories (see \cite{SGA}) in order to recover Theorem~\ref{1.11}
and Theorem~\ref{3.5}. Unfortunately this cannot be done in general because the
category of covering functors with finite fibre of $A$ may not be a
Galois category as explained in the following example:
\begin{example}
Let $A=kQ/I$ where $Q$ is equal to
\\
\null\hfill
\xymatrix{
 & 2 \ar@{->}[rd]^c & & 4 \ar@{->}[rd]^f &\\
1 \ar@{->}[ru]^b \ar@{->}[rr]_a & & 3 \ar@{->}[ru]^e \ar@{->}[rr]_d &&5
}
\hfill\null\\
and $I=< da,\ dcb +fea,\ fecb>$. Set $G=\mathbb{Z}/2=<\sigma|\sigma^2>$. Let $Q'$ be the quiver:\\
\null\hfill
\xymatrix{
& & 1 \ar@{->}[lldd]_b \ar@{->}[r]^a & \sigma 3 \ar@{->}[rdd]_{\sigma e}
\ar@{->}[rrdd]^{\sigma d}& & \\
&&\sigma 2 \ar@{->}[ru]_{\sigma c}&&& \\
2 \ar@{->}[rrdd]_{c}&  \sigma 1 \ar@{->}[ru]_{\sigma b}
\ar@{->}[ddr]^{\sigma a}&&& \sigma 4 \ar@{->}[ld]_{\sigma f} &5\\
&&&\sigma 5&&\\
&&3 \ar@{->}[ru]^{d} \ar@{->}[r]_{e}& 4\ar@{->}[rruu]_{f}&&
}
\hfill\null\\
and set:
\begin{equation}
I'=<\sigma\!d\:a,\ d\:\sigma\!a,\ dcb+\sigma\!f\:\sigma\!e\:a,\
\sigma\!d\:\sigma\!c\:\sigma\!b+fe\,\sigma\!a,\ fecb,\
\sigma\!f\:\sigma\!e\:\sigma\!c\:\sigma\!b>\notag
\end{equation}
Hence the natural mapping $p\colon (Q',I')\to (Q,I)$ ($x, \sigma
x\mapsto x$) is a Galois covering with group $G$. Therefore, if we set
 $A'=kQ'/I'$, then $p$ induces a Galois covering $F\colon
A'\to A$ with group $G$. If $u$ is a path in $Q'$ (resp. in $Q$) we
will write $\widetilde{u}$ (resp. $\widehat{u}$) for: $u\:mod\:I'\in A'$
 (resp. for: $u\:mod\:I\in A$).
Let us set $F'\colon A'\to A$ to be the
Galois covering with group $G$ as well and defined as follows:
\begin{enumerate}
\item[.] $F'(\widetilde{a})=F'(\widetilde{\sigma a})=\widehat{a}+\widehat{cb}$,
\item[.] $F'(\widetilde{x})=F'(\widetilde{\sigma x})=\widehat{x}$ for any
 $x\in\{b,c,d,e,f\}$.
\end{enumerate}
Assume that the category of the coverings of $A$ with finite fibre
is a Galois category. Hence this category admits finite products and the product of $F$ with $F'$ gives rise to a diagram:
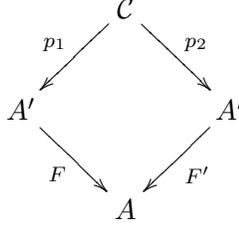
\null\hfill
\xymatrix{
&\c C\ar@{->}[rd]^{p_2} \ar@{->}[ld]_{p_1}&\\
 A'\ar@{->}[rd]_F && A'\ar@{->}[ld]^{F'}\\
&A &
}
\hfill\null\\
such that $F''=F\circ p_1=F'\circ p_2$ is a covering functor with fibre
 the product of the fibres of $F$ and $F'$. Hence $\c
 C_0=Q_0'\times_{Q_0}Q_0'=\bigcup_{x\in Q_0} \{(x,x),(x,\sigma x),(\sigma x,x),(\sigma
 x,\sigma x)\}$.
 Moreover, 
 Proposition~\ref{2.2} implies that $p_1$ and
 $p_2$ are covering functors as well.
Let us compute the lifting $u$ of $\widehat{a}\in\ _3A_1$
w.r.t. $F''$ and with source $(1,1)$. Using the lifting property of
 $p_1$ and $p_2$ we get:
\begin{enumerate}
\item[$\cdot$] $p_1(u_1)+p_1(u_2)=\widetilde{a}$ where $u_1+u_2\in\ _{(\sigma 3,3)}\c
 C_{(1,1)}\oplus\ _{(\sigma 3,\sigma 3)}\c C_{(1,1)}$,
\item[$\cdot$] $p_2(v_1)+p_2(v_2)=\widetilde{a}$ where $v_1+v_2\in\ _{(3,\sigma
 3)}\c C_{(1,1)}\oplus\ _{(\sigma 3,\sigma 3)}\c C_{(1,1)}$,
\item[$\cdot$] $p_2(v_3)+p_2(v_4)=\widetilde{cb}$ where $v_3+v_4\in\ _{(\sigma 3,3)}\c
 C_{(1,1)}\oplus\ _{(3, 3)}\c C_{(1,1)}$.
\end{enumerate}
Since $\widehat{a}=F(\widetilde{a})=F'(\widetilde{a}-\widetilde{cb})$, we infer
that $u=u_1+u_2=v_1+v_2-v_3-v_4$. Therefore $v_1=v_4=0$, $u_1=-v_3$ and
$u_2=v_2$.
Notice that $v_3\neq 0$ and $v_2\neq 0$ since $\widetilde{a}\neq 0$ and
$\widetilde{cb}\neq 0$. Therefore, the spaces $_{(\sigma 3,3)}\c
C_{(1,1)}$ and $_{(\sigma 3,\sigma 3)}\c C_{(1,1)}$ are non
zero. Since $p_1$ is a covering functor, we infer that  $p_1$ induces
an inclusion $_{(\sigma 3,\sigma 3)}\c C_{(1,1)}\oplus\ _{(\sigma
3,3)}\c C_{(1,1)}\hookrightarrow\ _{\sigma
 3}A'_1$ of a space of dimension at least $2$ in $_{\sigma
 3}A'_1=k.\widetilde{a}$. 
 This contradiction shows that the product of $F$ with $F'$ does not
 exist and that the category of coverings of $A$ with
 finite fibre is not necessarily a Galois category.
\end{example}
We end this study with a final remark concerning monomial
algebras. Recall that an algebra $A$ is monomial if it admits a
presentation $kQ/I_0\simeq A$ where $I_0$ is generated by a set of
paths. In such a case, $\pi_1(Q,I_0)\simeq\pi_1(Q)$ (the fundamental
group of $Q$) and therefore the fundamental group of any other
presentation of $A$ is a quotient of $\pi_1(Q,I_0)$. Thus,
Theorem~\ref{A} holds for $A$ monomial without hypothesis on the characteristic
of $k$ or on the double bypasses in $Q$. Hence we can wonder if
Theorem~\ref{B} holds for monomial algebras. This question will be
studied in a subsequent text. 
\section*{Acknowledgements}
The author would like to thank Claude Cibils, Mar\'ia Julia Redondo
and Andrea Solotar for carefully reading a previous version of this
text and for many helpful comments. He also acknowledges an anonymous referee,
Katie Coon and Guillaume Lamy for their help to improve the English of
the text.
\bibliographystyle{plain}
\bibliography{biblio}
\end{document}